\documentclass[12pt,reqno,twoside]{article}
\usepackage{amsmath,amsthm,amsfonts,amssymb,setspace,verbatim,cite,fancyhdr}

\usepackage[english]{babel}
\usepackage[utf8]{inputenc}
\usepackage{times}
\usepackage[T1]{fontenc}
\usepackage{mathrsfs}

\let\oldnsubseteq\nsubseteq 
\usepackage{mathabx} 
\let\nsubseteq\oldnsubseteq

\usepackage{enumitem}
\usepackage{accents}
\usepackage{fge}

\usepackage[top=15mm,bottom=13mm,left=22.8mm,right=22.8mm,includeheadfoot,headsep=8mm,footnotesep=5mm]{geometry} 

\usepackage[hypertexnames=false]{hyperref} 
\hypersetup{
pdfencoding=auto,
pageanchor=true,
plainpages=true,
pdftoolbar=true,
pdfmenubar=true,
pdfpagemode=UseOutlines,
pdffitwindow=false,
pdfstartview={FitH},
pdfnewwindow=true,
colorlinks=true,
linkcolor=black,
citecolor=black,
filecolor=black,
urlcolor=black,
linktoc=all,
breaklinks=true,
}

\usepackage{mdframed,xcolor}
\usepackage{xpatch}

\makeatletter
\xpatchcmd{\endmdframed}
  {\aftergroup\endmdf@trivlist\color@endgroup}
  {\endmdf@trivlist\color@endgroup\@doendpe}
  {}{}
\makeatother

\mdfdefinestyle{defStyle}
{linewidth=0.4mm,
skipabove=0.5\baselineskip,
skipbelow=0.9\baselineskip,
splittopskip=\topskip,
innertopmargin=0pt,
theoremspace={},
}

\newmdtheoremenv[linewidth=0.4mm,skipabove=0.9\baselineskip,skipbelow=0.9\baselineskip,splittopskip=\topskip,%
innertopmargin=0pt]{theorem}{Theorem}[section]
\newmdtheoremenv[linewidth=0.4mm,skipabove=0.5\baselineskip,skipbelow=0.9\baselineskip,splittopskip=\topskip,%
innertopmargin=0pt]{lemma}[theorem]{Lemma}
\newmdtheoremenv[linewidth=0.4mm,skipabove=0.5\baselineskip,skipbelow=0.9\baselineskip,splittopskip=\topskip,%
innertopmargin=0pt]{corollary}[theorem]{Corollary}
\newmdtheoremenv[linewidth=0.4mm,skipabove=0.5\baselineskip,skipbelow=0.9\baselineskip,splittopskip=\topskip,%
innertopmargin=0pt]{proposition}[theorem]{Proposition}
\mdtheorem[style=defStyle]{definition}[theorem]{Definition}
\newmdtheoremenv[style=defStyle]{hypothesis}[theorem]{Hypothesis}

\usepackage{thmtools}

\declaretheorem[style=definition,name=Remark,sibling=theorem]{remark}
\declaretheorem[style=definition,name=Example,qed=$\blacktriangle$,sibling=theorem]{example}

\numberwithin{equation}{section}

\makeatletter
\def\th@plain{%
  \thm@notefont{\bfseries}
  \itshape 
}

\def\th@definition{%
  \normalfont 
  \thm@notefont{\bfseries}
}
\makeatother

\renewcommand{\thefootnote}{\fnsymbol{footnote}}

\newcommand{\makepapertitle}{
\pdfbookmark[1]{Title page}{Title_page}
\thispagestyle{empty}\vspace*{8mm}
\begin{spacing}{1.1}
\LARGE\MakeUppercase{\papertitle}
\end{spacing}\vspace*{10mm}
}

\newcommand{\paperfirstauthor}{%
{\large\sc\firstauthor}\\[2mm] 
{\rm\firstaddress \\ 
\firstemail}\\[8mm]
}

\newcommand{\makeabstract}{%
\begin{minipage}{0.9\textwidth}
{\small {\sc Abstract.}
 \paperabstract
}
\end{minipage}\vfill
}

\newcommand{\MakeFirstPageOneAuthor}{
\begin{center}
  \makepapertitle
  \paperfirstauthor
  \makeabstract
  \begin{minipage}[t]{0.3\textwidth}
   \raggedleft {\bf Date (final version):} 
  \end{minipage}\hspace{0.01\textwidth}
  \begin{minipage}[t]{0.6\textwidth}
    \today\\ 
    (submitted on September 29, 2020;\\ accepted on February 4, 2021)
  \end{minipage}\\[4mm]
\noindent%
  \begin{minipage}[t]{0.3\textwidth}
    \raggedleft {\bf Running head:} 
  \end{minipage}\hspace{0.01\textwidth}
  \begin{minipage}[t]{0.6\textwidth}
    \runninghead
  \end{minipage}\\[4mm]
\noindent%
  \begin{minipage}[t]{0.3\textwidth}
    \raggedleft {\bf How to cite:} 
  \end{minipage}\hspace{0.01\textwidth}
  \begin{minipage}[t]{0.6\textwidth}
    J. Math. Anal. Appl. {\bf 499} (2021), no.~2, article~125054.\\[1mm]
   \footnotesize{\url{https://doi.org/10.1016/j.jmaa.2021.125054}}
  \end{minipage}\\
\bigskip
\vfill
\end{center}

\renewcommand{\thefootnote}{}
\footnotetext[2]{2010 {\it Mathematics Subject Classification:} \thesubjclass}
\footnotetext[3]{{\it Key words and phrases:} \thekeywords}
\setcounter{footnote}{0}
\renewcommand{\thefootnote}{{\bf\,\alph{footnote}\alph{footnote}\alph{footnote}\,}}

\setcounter{page}{0}
\newpage
}

\makeatletter
\newcommand{\rightorleftmark}{%
  \begingroup\protected@edef\x{\rightmark}%
  \ifx\x\@empty
    \endgroup\leftmark
  \else
    \endgroup\rightmark
  \fi}
\makeatother

\newcommand{\papertitle}%
{Eigenfunctions expansion for discrete symplectic systems with general linear dependence on spectral parameter}

\newcommand{\runninghead}%
{Eigenfunctions expansion for discrete symplectic systems}

\newcommand{\firstauthor}%
{Petr Zem{\'{a}}nek}
\newcommand{\firstauthorhead}%
{P. Zem{\'{a}}nek}

\newcommand{\firstaddress}
{Department of Mathematics and Statistics, Faculty of Science, Masaryk University \\
Kotl{\'{a}}{\v{r}}sk{\'{a}} 2, CZ-61137 Brno, Czech Republic}

\newcommand{\firstemail}%
{E-mail: zemanekp@math.muni.cz}

\newcommand{\paperabstract}%
{Eigenfunctions expansion for discrete symplectic systems on a finite discrete interval is established in the case of 
a general linear dependence on the spectral parameter as a significant generalization of the Expansion theorem given by 
Bohner, Do{\v{s}}l{\'{y}} and Kratz in [Trans. Amer. Math. Soc.~\textbf{361} (2009), 3109--3123]. Subsequently, 
an integral representation of the Weyl--Titchmarsh \texorpdfstring{$M(\la)$}{M(λ)}-function is derived explicitly by 
using a suitable spectral function and a possible extension to the half-line case is discussed. The main results 
are illustrated by several examples.}

\newcommand{\thekeywords}%
{Discrete symplectic system; eigenvalue; eigenfunction; expansion theorem; \texorpdfstring{$M(\la)$}{M(λ)}-function.}

\newcommand{\thesubjclass}%
{{\it Primary\/} 47B39; {\it Secondary\/} 39A12; 39A06; 34L10.}

\newcommand{\submittedto}%
{Journal of Mathematical Analysis and Applications}

\DeclareMathAccent{\wwtilde}{\mathord}{largesymbols}{"65}
\DeclareMathSymbol{\widetildesym}{\mathord}{largesymbols}{"65}
\newcommand\lowerwidetildesym{%
  \text{\smash{\raisebox{-1.3ex}{%
    $\widetildesym$}}}}
\newcommand\wtilde[1]{%
  \mathchoice
    {\accentset{\displaystyle\lowerwidetildesym}{#1}}
    {\accentset{\textstyle\lowerwidetildesym}{#1}}
    {\accentset{\scriptstyle\lowerwidetildesym}{#1}}
    {\accentset{\scriptscriptstyle\lowerwidetildesym}{#1}}
}

\newcommand\lowerwidetildesymW{%
  \text{\smash{\raisebox{-1.35ex}{$\widetildesym$}}}}
\newcommand\wtildeW[1]{%
    {\accentset{\scalebox{1.2}{\lowerwidetildesymW}}{#1}}
}

\DeclareMathSymbol{\widehatsym}{\mathord}{largesymbols}{"62}
\newcommand\lowerwidehatsym{%
  \text{\smash{\raisebox{-1.3ex}{%
    $\widehatsym$}}}}
\newcommand\what[1]{%
  \mathchoice
    {\accentset{\displaystyle\lowerwidehatsym}{#1}}
    {\accentset{\textstyle\lowerwidehatsym}{#1}}
    {\accentset{\scriptstyle\lowerwidehatsym}{#1}}
    {\accentset{\scriptscriptstyle\lowerwidehatsym}{#1}}
}

\DeclareMathAlphabet{\mthdtcl}{U}{dutchcal}{m}{n}
\DeclareMathAlphabet{\mathpzc}{OT1}{pzc}{m}{it}
\DeclareMathAlphabet{\msfsl}{U}{eus}{m}{n}

\renewcommand{\d}{\mathrm{d}}

\newcommand{\Ac}{\mathcal{A}}
\newcommand{\tAc}{\wtilde{\Ac}}
\newcommand{\Bc}{\mathcal{B}}
\newcommand{\tBc}{\wtilde{\Bc}}
\newcommand{\Cc}{\mathcal{C}}
\newcommand{\tCc}{\wtilde{\Cc}}
\newcommand{\Dc}{\mathcal{D}}
\newcommand{\tDc}{\wtilde{\Dc}}

\newcommand{\Ic}{\mathcal{I}}

\newcommand{\Jc}{\mathcal{J}}

\newcommand{\Oc}{\mathcal{O}}

\newcommand{\Sc}{\mathcal{S}}
\newcommand{\tSc}{\wtilde{\Sc}}
\newcommand{\hSc}{\what{\Sc}}

\newcommand{\Vc}{\mathcal{V}}
\newcommand{\tVc}{\wtilde{\Vc}}
\newcommand{\hVc}{\what{\Vc}}
\newcommand{\Wc}{\mathcal{W}}
\newcommand{\tWc}{\wtildeW{\Wc}}
\newcommand{\Xc}{\mathcal{X}}

\newcommand{\mL}{\mathscr{L}}

\newcommand{\Cbb}{\mathbb{C}}

\newcommand{\Nbb}{\mathbb{N}}
\newcommand{\Rbb}{\mathbb{R}}
\newcommand{\Sbb}{\mathbb{S}}
\newcommand{\tSbb}{\wtilde{\Sbb}}

\newcommand{\Zbb}{\mathbb{Z}}

\newcommand{\al}{\alpha}
\newcommand{\be}{\beta}

\newcommand{\la}{\lambda}
\newcommand{\hla}{\hat{\la}}

\newcommand{\bla}{\bar{\la}}

\newcommand{\de}{\delta}
\newcommand{\De}{\Delta}

\newcommand{\Ps}{\Psi}
\newcommand{\tPs}{\wtilde{\Ps}}
\newcommand{\hPs}{\what{\Ps}}
\newcommand{\Ph}{\Phi}

\newcommand{\Om}{\Omega}

\newcommand{\eps}{\varepsilon}
\newcommand{\Ga}{\Gamma}

\newcommand{\rh}{\rho}
\newcommand{\om}{\omega}
\newcommand{\si}{\sigma}

\newcommand{\stm}{\hspace*{0.2mm}\fgebackslash\hspace*{0.3mm}}

\usepackage{scalerel}
\usepackage{stackengine}

\stackMath

\newcommand{\hf}{\hat{f}}

\newcommand{\hu}{\hat{u}}

\newcommand{\hx}{\hat{x}}

\newcommand{\hy}{\hat{y}}

\newcommand{\tz}{\tilde{z}}

\newcommand{\tZ}{\wtilde{Z}}
\newcommand{\hz}{\hat{z}}
\newcommand{\hZ}{\what{Z}}

\newcommand{\ltp}{\ell^{\hspace{0.2mm}2}_{\Ps}}

\newcommand{\tltp}{\tilde{\ell}^{\hspace{0.3mm}2}_{\Ps}}
\newcommand{\tlt}[1]{\tilde{\ell}^{\hspace{0.3mm}2}_{#1}}

\newcommand{\sZbb}{{\scriptscriptstyle{\Zbb}}}
\newcommand{\Iz}{\Ic_\sZbb}

\newcommand{\Izp}{\Ic_\sZbb^+}

\newcommand{\oinftyZ}{[0,\infty)_\sZbb}
\newcommand{\onnZ}{[0,N+1)_\sZbb}
\newcommand{\onZ}{[0,N]_\sZbb}

\newcommand{\mmatrix}[1]{\begin{pmatrix} #1
  \end{pmatrix}}
\newcommand{\msmatrix}[1]{\left(\begin{smallmatrix} #1
  \end{smallmatrix}\right)}  

\newcommand{\qtextq}[1]{\quad\text{#1}\quad}
\newcommand{\qtext}[1]{\quad\text{#1 }\ }

\DeclareMathOperator{\re}{Re}
\DeclareMathOperator{\im}{Im}

\DeclareMathOperator{\ran}{ran}
\DeclareMathOperator{\rank}{rank}

\DeclareMathOperator{\tr}{tr}

\usepackage{mathtools,xparse}
\renewcommand{\.}{\hspace*{0.1 em}}
\DeclarePairedDelimiter\xnorm{\lVert}{\rVert}
\NewDocumentCommand{\norm}{som}
 {\IfBooleanTF{#1}
   {\xnorm*{#3}}
   {\IfNoValueTF{#2}
     {\mathopen{|\mkern-.8mu|}\.#3\.\mathclose{|\mkern-.8mu|}}
     {\xnorm[#2]{\.#3\.}}%
   }
 }
\DeclarePairedDelimiter\xinner{\langle}{\rangle}

\NewDocumentCommand{\xinnr}{som}
 {\IfBooleanTF{#1}
   {\xinner*{#3}}
   {\IfNoValueTF{#2}
     {\mathopen{\langle}\.#3\.\mathclose{\rangle}}
     {\xinner[#2]{\.#3\.}}%
   }
 }

\makeatletter
\def\inner{\@ifnextchar[{\@INNwith}{\@INNwithout}}
\def\@INNwith[#1]#2#3{\xinnr[#1]{#2,#3}}
\def\@INNwithout#1#2{\xinnr{#1,#2}}

\def\innerP{\@ifnextchar[{\@INNPwith}{\@INNPwithout}}
\def\@INNPwith[#1]#2#3{\xinnr[#1]{#2,#3}_\Ps}
\def\@INNPwithout#1#2{\xinnr{#1,#2}_\Ps}

\def\innerPN{\@ifnextchar[{\@INNPNwith}{\@INNPNwithout}}
\def\@INNPNwith[#1]#2#3{\xinnr[#1]{#2,#3}_{\Ps,N}}
\def\@INNPNwithout#1#2{\xinnr{#1,#2}_{\Ps,N}}

\def\abs{\@ifnextchar[{\@awith}{\@awithout}}
\def\@awith[#1]#2{{#1|}#2\.{#1|}}
\def\@awithout#1{|#1\.|}
 
\def\normS{\@ifnextchar[{\@Nwith}{\@Nwithout}}
\def\@Nwith[#1]#2{\norm[#1]{#2}_\si}
\def\@Nwithout#1{\norm{#1}_\si}

\def\normE{\@ifnextchar[{\@NEwith}{\@NEwithout}}
\def\@NEwith[#1]#2{\norm[#1]{#2}_2}
\def\@NEwithout#1{\norm{#1}_2}

\def\normA{\@ifnextchar[{\@NAwith}{\@NAwithout}}
\def\@NAwith[#1]#2{\norm[#1]{#2}_1}
\def\@NAwithout#1{\norm{#1}_1}

\def\normP{\@ifnextchar[{\@NPwith}{\@NPwithout}}
\def\@NPwith[#1]#2{\norm[#1]{#2}_{\Ps}}
\def\@NPwithout#1{\norm{#1}_{\Ps}}

\def\normW{\@ifnextchar[{\@NWwith}{\@NWwithout}}
\def\@NWwith[#1]#2{\norm[#1]{#2}_{\Wc}}
\def\@NWwithout#1{\norm{#1}_{\Wc}}

\def\normtW{\@ifnextchar[{\@NtWwith}{\@NtWwithout}}
\def\@NtWwith[#1]#2{\norm[#1]{#2}_{\scalebox{0.7}{$\tWc$}}}
\def\@NtWwithout#1{\norm{#1}_{\scalebox{0.6}{$\tWc$}}}

\newcommand{\Sla}[1]{\text{\rm(S$_{#1}$})}
\newcommand{\Slaf}[2]{\text{\rm(S$_{#1}^{#2}$)}}
\newcommand{\hSlaf}[2]{\text{\rm($\hat{\text{S}}_{#1}^{#2}$)}}

\DeclareMathOperator{\diag}{diag}
\DeclareMathOperator{\sgn}{sgn}

\DeclareMathOperator{\Eigen}{Eigen}
\DeclareMathOperator{\geom}{geom}
\DeclareMathOperator{\alg}{alg}

\newcommand{\ssc}[1]{{\scriptscriptstyle{#1}}}

\makeatletter
\newcommand{\ltxlabel}{\ltx@label}
\makeatother
\newcounter{GatherItemCounter}

\hypersetup{
pdfencoding=unicode,
pdftitle={\papertitle},
pdfauthor={\firstauthor},
pdfsubject={Department of Mathematics and Statistics, Faculty of Science, Masaryk University},
pdfkeywords={\thekeywords},
pdfsubject=\paperabstract,
pdfencoding=auto
}

\begin{document}


\MakeFirstPageOneAuthor



\section{Introduction}\label{S:intro}

Eigenfunctions expansion is a very classical topic in the spectral theory of linear operators on a~Hilbert 
space and the research in this field is still very active. Its origin can be traced back at least to works of 
Euler, d’Alembert, D. Bernoulli and, especially, of Fourier, see \cite{jbjF24,jbjF26,IGG72,jbjF09}. Subsequently, in 
the works \cite{cfS29,cfS36:i,cfS36:ii,cfS.jL37,jL37} of Sturm and Liouville, the theory of the eigenfunctions 
expansion was built in a more general fashion of regular boundary value problems associated with the second order linear 
differential equations. Its generalization to the case of a system of the first order differential equations was 
initiated by Hurwitz in \cite{waH21} and all of these results can be unified by using the linear Hamiltonian 
differential system~\eqref{E:cLHS} displayed below, see \cite{ccC22,aS21,gdB.reL23,gaB26,gaB38} and also
\cite{fvA64,rJ.rO.sN.cN.rF16,bcO69,laS99,wtR71:Wiley,wtR80}.

Although \cite{tA66,tA.sK72,vkD.hdD.rT02,sB01,cmB.lrM.kaM94} show that difference equations, with or without their 
related spectral theory, play also an important role in mathematical physics or continuum mechanics (see, for instance, 
\cite{tL.nP.gV19,tL.gV18,gV.jM13}), the literature covering discrete counterparts of expansion theorems discussed above 
seems to be, surprisingly, quite humbler. Some of them can be found, e.g., in \cite[Sections~4.4 and~6.8]{fvA64}, 
\cite[Section~2.5]{gT00}, \cite[Section~VII.1.10]{jmB68}, \cite[Section~2.5]{aJ95} or \cite{rdC22:BAMS,mhA.zsM.iaS12} 
for the second order Sturm--Liouville difference equations or three term recurrence relations (related to the Jacobi 
operator) and in a more general setting of the time scale calculus also in 
\cite{gsG07,gsG08,aH.eB09:NDST,faD.bpR07:JMAA,hT16}. But none of them deals with a~system. The only result of that type 
was published in~\cite{mB.oD.wK09} for a linear system with a very special dependence on the spectral parameter 
and it is recalled in Theorem~\ref{T:EE.admissible} below. Thus, in the present paper, we aim to extend the latter 
result to the case of general linear dependence on the spectral parameter.

More specifically, we develop an eigenfunctions expansion for a class of boundary value problems determined by the 
regular discrete symplectic system depending linearly on $\la\in\Cbb$, i.e.,
 \begin{equation*}\label{E:Sla}\tag{S$_\la$}
   z_{k}(\la)=(\Sc_k+\la\.\Vc_k)\,z_{k+1}(\la), \quad k\in\Iz,
 \end{equation*}
where $\Iz=\onZ\coloneq[0,N]\cap\Zbb$ is a finite discrete interval and, for every $k\in\Iz$, the coefficients 
$\Sc_k,\Vc_k$ are $2n\times2n$ matrices such that
 \begin{equation}\label{E:coeff.condition}
  \Sc_k^*\Jc\Sc_k=\Jc,\quad \Vc_k^*\Jc\Sc_k\ \text{ is Hermitian,} \qtextq{and} 
  \Vc_k^*\Jc\,\Vc_k=0
 \end{equation}
for $\Jc$ being the $2n\times2n$ orthogonal and skew-symmetric matrix 
$\Jc\coloneq\msmatrix{0 & -I\\ I & \phantom{-}0}$ with $n\times n$ blocks of zero and identity matrices. The first 
condition in~\eqref{E:coeff.condition} means that $\Sc_k$ are symplectic matrices and all conditions 
in~\eqref{E:coeff.condition} can be simultaneously written by using the matrix $\Sbb_k(\la)\coloneq \Sc_k+\la\.\Vc_k$ 
as the symplectic-type identity
 \begin{equation}\label{E:Sbb.symplectic}
  \Sbb_k^*(\bla)\.\Jc\.\Sbb_k(\la)=\Jc,
 \end{equation}
which is valid for all $k\in\Iz$ and $\la\in\Cbb$. This equality shows that $\Sbb_k(\la)$ is invertible with 
$\Sbb_k^{-1}(\la)=-\Jc\.\Sbb_k^*(\bla)\.\Jc$. Thus, \eqref{E:Sbb.symplectic} is equivalent to
$\Sbb_k(\bla)\.\Jc\.\Sbb_k^*(\la)=\Jc$ and similarly the conditions from~\eqref{E:coeff.condition} can be 
written with the superscript ``$*$'' applied to the third matrix instead of the first one. Consequently, solutions of 
systems~\eqref{E:Sla} and~\Sla{\nu} are connected via a~suitable Lagrange-type identity and, at the same time, 
system~\eqref{E:Sla} admits the equivalent ``operator-theoretic'' expression 
 \begin{equation}\label{E:equiv.Sla}
  \mL(z)_k\coloneq\Jc[z_k(\la)-\Sc_k\.z_{k+1}(\la)]
   =\la\.\Jc\.\Vc\.\.\Sbb_k^{-1}(\la)\.z_{k}(\la)=\la\.\.\Ps_k\.z_k(\la), \quad k\in\Iz,
 \end{equation}
where $\Ps_k\coloneq \Jc\.\Sc_k\.\Jc\.\Vc^*_k\.\Jc$ are $2n\times2n$ matrices such that
 \begin{equation}\label{E:Ps.assmpt}
   \Ps_k=\Ps_k^*\geq0 \qtextq{and} \Ps_k^*\,\Jc\,\Ps_k=0 \qtext{for all $k\in\Iz$.}
 \end{equation}
The additional assumption concerning the positive semidefiniteness of $\Ps_k$ is forced by the fact that it appears 
as the weight in the semi-inner product 
 \begin{equation*}
  \innerP{z}{u}\coloneq\sum_{k\in\Iz} z_k^*\,\Ps_k\,u_k \qtextq{and in the semi-norm}
  \normP{z}\coloneq\sqrt{\innerP{z}{z}},
 \end{equation*} 
which are naturally induced by system~\eqref{E:Sla}.

Besides the conditions for $\Vc_k$ or $\Ps_k$ displayed above, the linear dependence of~\eqref{E:Sla} on the parameter 
$\la$ is quite general. System~\eqref{E:Sla} is in the so-called ``time-reversed form'', because of the presence 
of the shift operator $k+1$ acting on the right-hand side instead of the left-hand side. Due to the absence of the 
shift operator in the semi-inner product $\innerP{\cdot}{\cdot}$, it was identified as the appropriate object for the 
study of the spectral theory of discrete symplectic systems, see \cite{slC.pZ15,pZ.slC16,pZ.slC:SAE2}. It can 
be easily transformed into the more traditional ``forward-time'' form 
 \begin{equation}\label{E:dss}
  z_{k+1}(\la)=\tSbb_k(\la)\,z_k(\la), \quad k\in\Iz,
 \end{equation}
with $\tSbb_k(\la)=\tSc_k+\la\tVc_k\coloneq\Sbb^{-1}_k(\la)$. System~\eqref{E:dss} was introduced 
in~\cite{rSH.pZ14:JDEA} and after minor modifications we can adapt all known results to the setting of~\eqref{E:Sla}. 
The very special case of~\eqref{E:dss} with $\tPs_k=\diag\{\tWc_k,0\}$ for $\tWc_k=\tWc_k^*\geq0$ being $n\times n$ 
matrices on $\Iz$ was intensively studied in the last decade as the natural extension of any even-order 
Sturm--Liouville difference equation, see e.g.,
\cite{slC.pZ10,mB.oD.wK09,oD.wK07:JDEA,mB.sS10,oD09,oD.wK10:JDEA,jvE10:AML}. In addition, discrete symplectic systems 
are also known as the proper discrete analogue of the linear Hamiltonian differential systems related, e.g., by the 
symplectic structure of their fundamental matrices. Their systematic study was initiated in~\cite{mB.oD97,cdA.acP96} 
and they naturally arise, e.g., in the discrete calculus of variations, numerical integration schemes or in the theory 
of continued fractions. For a~thorough treatise on these systems we refer the reader to the monograph 
\cite{oD.jE.rsH19}.

If we partition the matrix $\tSc_k=\msmatrix{\tAc_k & \tBc_k\\ \tCc_k & \tDc_k}$ into $n\times n$ blocks and similarly 
put $z_k=(x_k^*,u_k^*)^*$, then the special diagonal structure of $\Ps_k$ mentioned above enables us to write the 
system as the pair of 
equations
 \begin{equation}\label{E:dss.blocks}
  \left.
   \begin{aligned}
    x_{k+1}(\la)&=\tAc_k\.x_k(\la)+\tBc_k\.u_k(\la),\\
    u_{k+1}(\la)&=\tCc_k\.x_k(\la)+\tDc_k\.u_k(\la)-\la\.\tWc_k\.x_{k+1}(\la),
   \end{aligned}
  \,\right\}
 \end{equation}
for which the parameter $\la$ appears only in the second-half of the system. It means that the space of all 
{\it admissible sequences} for~\eqref{E:dss.blocks}, i.e., sequences satisfying only the first equation 
of~\eqref{E:dss.blocks}, is constant with respect to $\la$. This observation was crucial in the proof of the 
{\it Rayleigh principle} for system~\eqref{E:dss.blocks} with real-valued coefficient matrices, from which the 
eigenfunctions expansion can be easily deduced for {\it finite eigenfunctions} of~\eqref{E:dss.blocks} being 
$\Wc$-nontrivial solutions of system~\eqref{E:dss.blocks} with zero endpoints, i.e., they satisfy~\eqref{E:dss.blocks} 
for some $\la\in\Cbb$ and their first components are such that
 \begin{equation}\label{E:eigen.blocks}
  x_0(\la)=0=x_{N+1}(\la), \qtextq{and} \Wc_k\.x_{k+1}(\la)\not\equiv0 \text{ on $\onZ$}.
 \end{equation}
The latter result reads as follows and it can be found in~\cite[Proposition~2]{oD.wK07:JDEA} 
and~\cite[Theorem~4.7]{mB.oD.wK09}, see also \cite[Section~4]{oD.wK10:JDEA}.

\begin{theorem}\label{T:EE.admissible}
 Let $\Iz=\onZ$ be a finite discrete interval and real-valued matrices $\tSc_k\in\Rbb^{2n\times2n}$ and 
 $\tWc_k\in\Rbb^{n\times n}$ be such that $\tSc_k^*\.\Jc\.\tSc_k=\Jc$ and $\tWc_k=\tWc_k^*\geq0$ for all $k\in\Iz$. 
 Then the number of all  linearly independent finite eigenfunctions of system~\eqref{E:dss.blocks}, denoted by $r$, is 
 finite and satisfies the estimate
  \begin{equation*}
   0\leq r\leq\sum_{k=0}^{N-1} \rank\tWc_k\leq nN.
  \end{equation*}
 Furthermore, these eigenfunctions form a complete orthonormal system in the linear space of all admissible sequences 
 of~\eqref{E:dss.blocks} with zero endpoints, i.e., the first component of any admissible sequence 
 $\{\hz_k\}_{k=0}^{N+1}$ satisfying $\hx_0=0=\hx_{N+1}$ can be ``expressed'' by using the first components of linearly 
 independent finite eigenfunctions $z^{[1]},\dots,z^{[r]}$ as the sum $\sum_{j=1}^r c_j\.x^{[j]}$ with the coefficients 
 $c_j\coloneq \sum_{k\in\Iz} x^{[j]*}_{k+1}\.\tWc_k\.\hx_{k+1}$ for all $j=1,\dots,r$ in the sense
  \begin{equation*}
   \normtW[\Big]{\hx-\sum_{j=1}^r c_j\.x_j}=0,
   \qtextq{where the semi-norm is given by} \normtW{x}\coloneq \sqrt{\sum_{k\in\Iz} x_{k+1}^*\.\tWc_k\.x_{k+1}}.
  \end{equation*}
 In particular, if $\tWc_k\equiv I$, then we have even the expansion
  \begin{equation*}
   \hx_k=\sum_{j=1}^r c_j\.x^{[j]}_k \qtext{for all $k\in[0,N+1]_\sZbb$ and } 
   c_j=\sum_{k\in\Iz} x^{[j]*}_{k+1}\.\hx_{k+1}.
 \end{equation*}
\end{theorem}

Theorem~\ref{T:EE.admissible} is a discrete counterpart of the expansion theorems given in \cite[Theorem~1]{wK95:JLMS} 
and \cite[Theorem~4.3]{wK.rSH12} for the linear Hamiltonian (or canonical) differential system
 \begin{equation}\label{E:cLHS}\tag{H$_\la$}
  -\Jc z'(t,\la)=[H(t)+\la\.W(t)]\.z(t,\la),\quad t\in[a,b],
 \end{equation}
where $H(t)$ and $W(t)\geq0$ are $2n\times 2n$ Hermitian matrices for all $t\in[a,b]$ with $W(t)$ possessing 
a~similar block-diagonal structure as in the case of~\eqref{E:dss.blocks}. These results were unified and their 
range was further extended by using the time scale calculus in \cite[Theorem~4.7]{rSH.vZ11}. However, for general 
linear dependence on $\la$ the space of admissible sequences is no longer independent of $\la$, 
which significantly complicates any generalization of the Rayleigh principle and, consequently, the expansion theorem 
to that case. Therefore, in the present paper we do not follow the approach based on the Rayleigh 
principle. Instead of that we aim to establish a discrete analogue of the expansion theorem presented 
in~\cite[Chapter~9]{fvA64} for system~\eqref{E:cLHS} with an arbitrary positive semidefinite matrix-valued function 
$W(\cdot)$. More precisely, we utilize the Weyl--Titchmarsh theory from~\cite{rSH.pZ14:JDEA} to expand any solution of 
the nonhomogeneous boundary value problem
 \begin{equation}\label{E:EE15a.intro}
  \mL(z)_k=\Ps_k\.f_k,\quad k\in\Iz,\quad \al\.z_0=0=\be\.z_{N+1}
 \end{equation}
by using {\it eigenfunctions} of~\eqref{E:Sla}, i.e., nontrivial solutions of the eigenvalue problem
 \begin{equation}\label{E:EP}
  \eqref{E:Sla},\quad \la\in\Cbb,\quad k\in\Iz,\quad \quad \al\.z_0(\la)=0=\be\.z_{N+1}(\la)
 \end{equation}
with the boundary conditions in~\eqref{E:EE15a.intro} and~\eqref{E:EP} being the same and determined by an 
arbitrary pair of matrices $\al,\be\in\Ga$ from the set
 \begin{equation}\label{E:Ga.def}
  \Ga\coloneq\{\al\in\Cbb^{n\times2n}\mid \al\.\al^*=I,\ \al\.\Jc\al^*=0\}.
 \end{equation} 
This alternative of Theorem~\ref{T:EE.admissible} can be stated as follows and it is established in 
Section~\ref{S:main}, see Theorem~\ref{T:main.expansion} for more precise formulation and also 
Corollary~\ref{C:expansion}. Albeit instead of finite eigenfunctions it formally deals with the ``standard'' notion of 
eigenfunctions introduced above, it relies on the so-called {\it Weak Atkinson condition}, which guarantees that every 
eigenfunction satisfies $\Ps_k z^{[j]}_k\not\equiv0$ on $\Iz$, i.e., they are finite in the sense 
of~\eqref{E:eigen.blocks} and all of them truly take place in the expansion. 

\begin{theorem}\label{T:main.intro}
 Let $\Iz=\onZ$ be a finite discrete interval, matrices $\al,\be\in\Ga$ be given, and $\Sc_k$ and $\Vc_k$ be 
 $2n\times2n$ complex-valued matrices such that all conditions in~\eqref{E:coeff.condition} and $\Ps_k\geq0$ hold for 
 all $k\in\Iz$. Then the number of all linearly independent eigenfunctions of~\eqref{E:Sla}, denoted by $r$, is finite 
 and satisfies the estimate
  \begin{equation*}
   0\leq r\leq 2n(N+2).
  \end{equation*}
 If, in addition, the Weak Atkinson condition is satisfied, then we have yet another upper bound 
 $r\leq\min\big\{n(N+1),\sum_{k\in\Iz}\rank\Ps_k\big\}$. Furthermore, for any sequence $\{f_k\}_{k=0}^{N+1}$ any 
 solution $\{\hz_k\}_{k=0}^{N+1}$ of problem~\eqref{E:EE15a.intro} admits the representation as the finite sum
 $\sum_{j=1}^r c_j\.z^{[j]}$ in the sense
  \begin{equation*}
   \normP[\Big]{\hz-\sum_{j=1}^r c_j\.z^{[j]}}=0,
  \end{equation*}
 where $c_j\coloneq \sum_{k\in\Iz} z^{[j]*}_k\.\Ps_k\.\hz_k$ for all $j=1,\dots,r$ with $z^{[1]},\dots,z^{[r]}$ being 
 linearly independent eigenfunctions of system~\eqref{E:Sla}.
\end{theorem}
A~relation between Theorems~\ref{T:EE.admissible} and~\ref{T:main.intro} is discussed in detail later in 
Remark~\ref{R:EE.theorems.compare}. Furthermore, in our treatise we utilize some properties of the associated 
Weyl--Titchmarsh $M(\la)$-function and, as a~``by product'' of Theorem~\ref{T:main.intro}, we derive its explicit 
integral representation in Theorem~\ref{T:Im.Mla.integral}, which is based on the Riemann--Stieltjes integral with 
respect to a suitable spectral function. Finally, let us note that the basic notation and assumptions used throughout 
the manuscript are summarized in the following Section~\ref{S:prelim}, which contains also several preliminary results.
\medskip

\section{Preliminaries}\label{S:prelim}

Throughout the paper, we employ the following notation and conventions. The real and imaginary parts of any 
$\la\in\Cbb$ are, respectively, denoted by $\re\la$ and $\im\la$, i.e., $\re\la\coloneq (\la+\bla)/2$ and 
$\im\la\coloneq (\la-\bla)/(2i)$. The symbols $\Cbb_+$ and $\Cbb_-$ mean, respectively, the upper and lower complex 
planes, i.e., $\Cbb_+\coloneq\{\la\in\Cbb\mid \im(\la)>0\}$ and $\Cbb_-\coloneq\{\la\in\Cbb\mid \im(\la)<0\}$. All 
matrices are considered over the field of complex numbers $\Cbb$. For $r,s\in\Nbb$ we denote by $\Cbb^{r\times s}$ the 
space of all complex-valued $r\times s$ matrices and $\Cbb^{r\times 1}$ will be abbreviated as $\Cbb^r$. In particular, 
the $r\times r$ \emph{identity} and \emph{zero matrices} are written as $I_r$ and $0_r$, where the subscript is 
omitted whenever it is not misleading (for simplicity, the zero vector is also written as $0$). For a given 
matrix $M\in\Cbb^{r\times s}$ we indicate by $M^*$, $\ker M$, $\ran M$, $\det M$, $\rank M$, $\re M\coloneq (M+M^*)/2$, 
$\im M\coloneq (M-M^*)/(2i)$, $M\geq0$, and $M>0$ respectively, its conjugate transpose, kernel, range (or image), 
determinant, rank, Hermitian components (real and imaginary parts), positive semidefiniteness, and positive 
definiteness. Finally, the symbol $\diag\{M_1,\dots,M_n\}$ stands for the block-diagonal matrix with matrices 
$M_1,\dots,M_n$ located on the main diagonal.

If $\Ic$ is an open or closed interval in $\Rbb$, then $\Iz\coloneq \Ic\cap\Zbb$ denotes the corresponding discrete 
interval. In particular, with $N\in\Nbb\cup\{0,\infty\}$, we shall be interested in discrete intervals of the form 
$\Iz\coloneq [0,N+1)_\Zbb$, in which case we define $\Izp\coloneq [0,N+1]_\sZbb$ with the understanding that 
$\Izp=\Iz$ when $N=\infty$, i.e., we deal only with discrete intervals $\Iz$ which are finite or unbounded above. 
By $\Cbb(\Iz)^{r\times s}$ we denote the space of sequences defined on $\Iz$ of complex $r\times s$ matrices, where 
typically $r\in\{n,2n\}$ and $1\leq s\leq2n$. In particular, we write only $\Cbb(\Iz)^r$ in the case $s=1$. If 
$M\in\Cbb(\Iz)^{r\times s}$, then $M(k)\coloneq M_k$ for $k\in\Iz$ and if $M(\la)\in\Cbb(\Iz)^{r\times s}$, then 
$M(\la,k)\coloneq M_k(\la)$ for $k\in\Iz$ with $M_k^*(\la)\coloneq[M_k(\la)]^*$. If $M\in\Cbb(\Iz)^{r\times s}$ and 
$L\in\Cbb(\Iz)^{s\times p}$, then $MN\in\Cbb(\Iz)^{r\times p}$, where $(MN)_k\coloneq M_k N_k$ for $k\in\Iz$. Finally, 
the forward difference operator acting on $\Cbb(\Iz)^{r\times s}$ is denoted by $\De$ with $(\De z)_k\coloneq\De z_k$ 
and we also put $z_k \big|_{m}^n\coloneq z_n - z_m$.

In the final part of this paper we utilize the Riemann--Stieltjes integral introduced by Stieltjes in 1894, see e.g., 
\cite{thH38}. In particular, we employ the following basic proposition, which is a direct consequence of the 
definition of this type of integral.

\begin{proposition}\label{P:RS-integral}
 Let $f:[a,b]\to\Rbb$ be a given function and $\tau:[a,b]\to\Rbb$ be a step function with discontinuities at 
 $a\leq x_1<x_2<\dots<x_n\leq b$. If, for every $k\in\{1,2,\dots,n\}$, either $f$ or $\tau$ is continuous from the 
 left at $x_k$ and, simultaneously, either $f$ or $\tau$ is continuous from the right at $x_k$, then $f$ is 
 Riemann--Stieltjes integrable on $[a,b]$ with respect to $\tau$ and it holds
  \begin{equation*}
   \int_a^b f(x)\.\d\tau(x)=\sum_{k=1}^n f(x_k)\.[\tau(x_k^+)-\tau(x_k^-)],
  \end{equation*}
 where $\tau(x_k^\pm)$ stands for the corresponding one-sided limits with $\tau(a^-)\coloneq\tau(a)$ and 
 $\tau(b^+)\coloneq\tau(b)$ whenever $x_1=a$ and $x_n=b$, respectively.
\end{proposition}

Let us recall that the primary object of our present treatise is represented by the time-reversed discrete symplectic 
system~\eqref{E:Sla} with the coefficient matrices satisfying~\eqref{E:coeff.condition} or \eqref{E:Sbb.symplectic} 
and/or \eqref{E:Ps.assmpt}. Since it can be equivalently written as~\eqref{E:equiv.Sla} and these two systems are 
mutually connected via $\Ps_k=\Jc\.\Sc_k\.\Jc\.\Vc^*_k\.\Jc$ or $\Vc_k=-\Jc\.\Ps_k\.\Sc_k$, it suffices to deal
only with one of the pairs $\{\Sc_k,\Vc_k\}$ and $\{\Sc_k,\Ps_k\}$ as it is done in the following hypothesis, which is 
tacitly assumed to hold whenever we speak about system~\eqref{E:Sla} or its nonhomogeneous counterpart
 \begin{equation}\label{E:Slaf}\tag{S$_\la^f$}
  z_k(\la)=\Sbb_k(\la)\,z_{k+1}(\la)-\Jc\,\Ps_k\,f_k, \quad k\in\Iz,
 \end{equation}
if not stated otherwise as in Theorem~\ref{T:hilbert}. It should be also pointed out that the existence of a 
unique solution $z(\la)\in\Cbb(\Izp)^{2n}$ of any initial value problem associated with system~\eqref{E:Sla} is easily 
guaranteed by equality~\eqref{E:Sbb.symplectic}, which yields $\abs{\det\Sbb_k(\la)}=1$ for all $\la\in\Cbb$ and 
$k\in\Iz$.
 
\begin{hypothesis}\label{H:basic}
 A number $n\in\Nbb$ and a finite discrete interval $\Iz=\onZ$ are given and we have $2n\times2n$ matrix-valued 
 sequences $\Sc,\Ps\in\Cbb(\Iz)^{2n\times2n}$ such that
  \begin{equation}\label{E:H.basic}
   \Sc_k^*\Jc\Sc_k=\Jc,\quad \Ps_k^*=\Ps_k, \quad \Ps_k^*\,\Jc\,\Ps_k=0, \qtextq{and} \Ps_k\geq0 
   \qtext{for all $k\in\Iz$.}
  \end{equation}
 Moreover, we define $\Sbb_{k}(\la)\coloneq\Sc_{k}+\la\Vc_{k}$ with $\Vc_k\coloneq -\Jc\,\Ps_k\,\Sc_k$ for all 
 $k\in\Iz$.
\end{hypothesis}

We emphasize that the singularity of matrices $\Ps_k$ contained in~\eqref{E:H.basic}(iii) is truly necessary and cannot 
be avoided as in the continuous-time case. It is a consequence of the possible expression of~\eqref{E:Sla} as 
a~linear perturbation of a certain ``operator'' given in~\eqref{E:equiv.Sla} and the Lagrange-type identity being 
``linear'' with respect to the parameter $\la$, see Theorem~\ref{T:Lagrange} below. Consequently, the set $\Cbb(\Izp)$ 
together with $\innerP{\cdot}{\cdot}$ is never an inner product and to get a~Hilbert space associated 
with~\eqref{E:Sla} 
we need to consider equivalences classes of sequences, compare with \cite[Lemma~2.5]{yS06}.

\begin{theorem}\label{T:hilbert}
 Let $\Iz=\onnZ$ be a finite or infinite discrete interval. Then
  \begin{equation*}
   \ltp=\ltp(\Iz)\coloneq\Big\{z\in\Cbb(\Izp)^{2n}\mid \normP{z}\coloneq \sqrt{\innerP{z}{z}}<\infty  \Big\}
  \end{equation*}
 is always only a semi-inner product space with respect to $\innerP{\cdot}{\cdot}$ and $\dim\ltp=2n(N+2)$, while the 
 corresponding space of equivalence classes
  \begin{equation*}
   \tltp=\tltp(\Iz)\coloneq \ltp\big/\big\{z\in\Cbb(\Izp)^{2n}\mid \ \normP{z}=0\big\}
  \end{equation*}
 is a Hilbert space with $\dim\tltp=\sum_{k\in\Iz} \rank \Ps_k$ and the inner product 
 $\innerP{[z]}{[y]}\coloneq \innerP{z}{y}$, which does not depend on the choice of representatives $z\in[z]\in\tltp$ 
 and $y\in[y]\in\tltp$. In particular, in the case  $\Iz=\oinftyZ$ we have $\dim\ltp=\infty$, while $\dim\tltp=\infty$ 
 holds if and only if, in addition, $\Ps_k\neq0$ for infinitely many $k\in\Iz$.
\end{theorem}
\begin{proof}
 The first and the last parts are quite obvious as well as the fact that the space $\tltp$ is an inner product space. 
 The independence of $\innerP{[\cdot]}{[\cdot]}$ follows from the fact that $z^{[1]},z^{[2]}\in[z]\in\tltp$ if and only 
 if $z^{[1]}-z^{[2]}\in[0]$ or equivalently $\Ps_k^{1/2}\.z_k^{[1]}=\Ps_k^{1/2}\.z_k^{[2]}$ for all $k\in\Iz$. It 
 remains to prove the completeness of the space $\tltp$ and its dimension. Let us start with the simplest case
 $\Ps_k=\diag\{\hPs_k,0\}$ for all $k\in\Iz$ with $\hPs_k>0$ being an $r_k\times r_k$ matrix for some 
 $r_k\in\{0,\dots,2n\}$. If we take an arbitrary Cauchy sequence $\{[z^{(p)}]\}_{p\in\Nbb}$ of classes 
 $[z^{(1)}],[z^{(2)}],\dots\in\tltp$, i.e., for any $\eps>0$ there is a number $M\in\Nbb$ such that for all $m,j\geq M$ 
 we have
  \begin{equation*}
   \normP[\big]{[z^{(m)}]-[z^{(j)}]}<\eps,
  \end{equation*}
 then the sequences $y^{(p)}\coloneq \Ps^{1/2} z^{(p)}\in\Cbb(\Izp)^{2n}$, $p\in\Nbb$, each of them being determined by 
 arbitrary representatives $z^{(p)}\in[z^{(p)}]$, obviously satisfy
  \begin{equation*}
   \norm[\big]{y^{(p)}}_{I}=\sum_{k\in\Iz} y^{(p)*}\.y^{(p)}=\normP[\big]{z^{(p)}}.
  \end{equation*}
 Consequently,
  \begin{equation*}
   \norm[\big]{y^{(m)}-y^{(j)}}_I=\normP[\big]{z^{(m)}-z^{(j)}},
  \end{equation*}
 which means that the corresponding equivalence classes $[y^{(1)}],[y^{(2)}],\dots$ belong to the space $\tlt{I}$ and 
 they form a~Cauchy sequence $\{[y^{(p)}]\}_{p\in\Nbb}$ in $\tlt{I}$. If $\Iz=\oinftyZ$, then $\tlt{I}$ is the Hilbert 
 space of square summable sequences on $\Iz$, while in the finite case it can be seen as the Hilbert space of sequences 
 defined on $\Iz$, each of them gives rise to a class of sequences, which differ only at $N+1$. In any case, there is 
 a unique $[\hy]\in\tlt{I}$ such that $\norm{[y^{(p)}]-[\hy]}_I\to0$ as $p\to\infty$. If we take an arbitrary 
 representative $\hy\in[\hy]$ and define
  \begin{equation*}
   \hz_k\coloneq \diag\{\hPs_k^{-1/2},0\}\,\hy_k,\quad k\in\Izp,
  \end{equation*}
 then the corresponding class $[\hz]$ belongs to $\tltp$ and the equality
 $\normP{[z^{(p)}]-[\hz]}=\norm{[y^{(p)}]-[\hy]}_I$ yields the convergence $[z^{(p)}]\to[\hz]$ in $\tltp$ as 
 $p\to\infty$. This proves the completeness of the space $\tltp$ in the given special case. 
 
 In a general case, the positive semidefiniteness of $\Ps_k$ implies the existence of a $2n\times2n$ unitary matrix 
 $\Om_k$ such that
  \begin{equation*}
   \Om^*_k\.\Ps_k\.\Om_k=\diag\{\hPs_k,0\}, \quad k\in\Iz,
  \end{equation*}
 where $\hPs_k>0$ is a suitable $r_k\times r_k$ matrix with $r_k\coloneq\rank\Ps_k$. If $\{[z^{(p)}]\}_{p\in\Nbb}$ is 
 again a~Cauchy sequence in $\tltp$ and $\{[u^{(p)}]\}_{p\in\Nbb}$ is the associated sequence of classes determined by 
 $u^{(p)}\coloneq \Om^*\.z^{(p)}$, then $\norm{[u^{(p)}]-[u^{(q)}]}_{\Om^*\Ps\Om}=\normP{[z^{(p)}]-[z^{(q)}]}$
 implies that $[u^{(p)}]\in\tlt{\Om^*\Ps\Om}$. Therefore, $\{[u^{p}]\}_{p\in\Nbb}$ is a~Cauchy sequence 
 in $\tlt{\Om^*\Ps\Om}$, which together with the previous part yields its convergence as well as the 
 convergence of $\{[z^{(p)}]\}_{p\in\Nbb}$ in $\tltp$ by using the inverse transformation, i.e., the inner product 
 space $\tltp$ is complete for any $\Ps\in\Cbb(\Iz)^{2n\times2n}$ satisfying the basic assumptions given 
 in~\eqref{E:H.basic}.
 
 Finally, we calculate the dimension of $\tltp$. First of all, we note that $[z]=[0]$ if and only if 
 $\Ps_k\.z_k\equiv0$ on $\Iz$ for all representatives $z\in[z]$, i.e., by using the previous notation, 
  \begin{equation*}
   \Om_k\.\diag\{\hPs_k,0\}\.\Om_k^*\.z_k=\Ps_k\.z_k=0\qtextq{or only} \diag\{\hPs_k,0\}\.\Om_k^*\.z_k=0 
   \qtext{for all $k\in\Iz$.}
  \end{equation*}
 Since for every $k\in\Iz$ there are $r_k$ linearly independent vectors $y_k\in\Cbb^{2n}$ such that 
 $\diag\{\hPs_k,0\}\.y_k\neq0$, the unitary property of $\Om_k$ implies the existence of $r_k$ linearly independent 
 vectors $z_k\coloneq \Om_k\.y_k$ satisfying $\Ps_k\.z_k\neq0$. Thus, we have precisely $r_k$ ``fundamental sequences'' 
 of the form $\{0,\dots,0,z_k,0,\dots\}$, each of which gives rise to a~distinct equivalence class, i.e., it holds 
 $\dim\tltp=\sum_{k\in\Iz} r_k$.
\end{proof}

The following identity, mentioned in the introduction, is one of the most important tool in the whole theory related 
to system~\eqref{E:Sla}. For its proof we refer to~\cite[Theorem~2.5]{slC.pZ15}.

\begin{theorem}[Extended Lagrange formula]\label{T:Lagrange}
 Let $\Iz$ be a finite or infinite discrete interval, numbers $\la,\nu\in\Cbb$, $m\in\{1,\dots,2n\}$ be given 
 as well as sequences $f,g\in\Cbb(\Iz)^{2n\times m}$. If $z(\la),u(\nu)\in\Cbb(\Izp)^{2n\times m}$ are solutions of 
systems~\eqref{E:Slaf} and \Slaf{\nu}{g}, respectively, then for 
 any $k,s,t\in\Iz$ such that $s\leq t$, we have 
  \begin{align}
   \De[z_k^*(\la)\,\Jc u_k(\nu)]&=(\bla-\nu)\,z_k^*(\la)\,\Ps_k\,u_k(\nu)+f_k^*\,\Ps_k\,u_k(\nu)
                                                            -z_k^*(\la)\,\Ps_k\,g_k,\notag\\
   &\hspace*{-33mm} z_k^*(\la)\,\Jc u_k(\nu)\big|_{s}^{t+1}
     =\sum_{k=s}^{t}\big\{(\bla-\nu)\,z_k^*(\la)\,\Ps_k\,u_k(\nu)+f_k^*\,\Ps_k\,u_k(\nu)
                                                            -z_k^*(\la)\,\Ps_k\,g_k\big\}.\label{E:lagrange}
  \end{align}
 Especially, if $\nu=\bla$ and $f\equiv0\equiv g$, then we get the Wronskian-type identity
  \begin{equation}\label{E:wronski.id}
   z_k^*(\la)\,\Jc\,u_k(\bla)=z_0^*(\la)\,\Jc\,u_0(\bla),\quad k\in\Izp.
  \end{equation}
\end{theorem} 

Consequently, if $\la\in\Cbb$ is arbitrary and $\Ph(\la)\in\Cbb(\Izp)^{2n\times2n}$ is a fundamental matrix of 
system~\eqref{E:Sla} such that it possesses the symplectic-type property $\Ph_s^*(\bla)\.\Jc\.\Ph_s(\la)=\Jc$ at some 
$s\in\Izp$, then equality~\eqref{E:wronski.id} implies that
 \begin{equation*}\label{E:Ph.symlectic}
  \Ph_k^*(\la)\.\Jc\.\Ph_k(\bla)=\Jc \qtextq{and}
  \Ph_k^{-1}(\la)=-\Jc\.\Ph_k^*(\bla)\.\Jc  \qtext{for all $k\in\Izp$.}
 \end{equation*}
In particular, it is satisfied for the fundamental matrix $\Ph(\la)=\big(\hZ(\la),\tZ(\la)\big)$ composing of 
the pair of linearly independent solutions $\hZ,\tZ\in\Cbb(\Izp)^{2n\times n}$ determined by the initial conditions
 \begin{equation}\label{E:Ph.via.al}
  \hZ_0(\la)=\al^* \qtextq{and} \tZ_0(\la)=-\Jc\al^*
 \end{equation}
for any $\al\in\Ga$ from the set defined in~\eqref{E:Ga.def}. It should be also noted that the symplectic property of 
the matrix $\big(\al^*,\,-\Jc\al^*\big)$ yields the inversion formula $\big(\al^*,\,-\Jc\al^*\big)^{-1}=\msmatrix{\al\\ 
\al\.\Jc}$ and $\ker\al=\ran\Jc\al^*$. The latter equality shows that the second half of this fundamental matrix shall 
be very useful in the decision whether $\la$ is an eigenvalue of~\eqref{E:Sla} or not, i.e., if there is a nontrivial 
solution of~\eqref{E:EP}. But before a~deeper analysis of this observation we need to introduce the following {\it Weak 
Atkinson condition}, which was identified in~\cite{rSH.pZ14:JDEA} as the minimal necessary assumption for any 
reasonable 
development of the Weyl--Titchmarsh theory for discrete symplectic systems or (in the continuous-time setting) for 
linear Hamiltonian differential system, see also~\cite{rSH.pZ14:AMC}. For the purpose of the paper at hand, this 
condition is still absolutely sufficient and we do not need to follow the traditional approach based on a similar 
assumption concerning all nontrivial solutions (the so-called {\it Strong Atkinson condition} or {\it definiteness 
condition}) as it can be found, e.g., in~\cite[Inequality~(9.1.6)]{fvA64}. For completeness, let us emphasize that in 
the rest of this section $\Iz$ stands for a~finite discrete interval only.

\begin{hypothesis}[Weak Atkinson condition]\label{H:WAC}
 For a given $\al\in\Ga$ there exists $\la\in\Cbb$ such that every nontrivial linear combination 
 $z(\la)\in\Cbb(\Izp)^{2n}$ of the columns of $\tZ(\la)$ satisfies
  \begin{equation}\label{E:Atk.ineq}
   \sum_{k=0}^{N} z_k^*(\la)\.\Ps_k\.z_k(\la)>0.
  \end{equation}
\end{hypothesis}

Although it seems, at first sight, that the present formulation of the Weak Atkinson condition is quite different than 
in~\cite[Hypothesis~2.7]{rSH.pZ14:JDEA}, we show quite easily that they are equivalent, i.e., the required property is 
either true for all $\la\in\Cbb$ or for none and eigenfunctions are only multiples of $\tZ(\la)$.

\begin{lemma}\label{L:definite}
 Let $\al\in\Ga$ be given and $\la\in\Cbb$ be such that Hypothesis~\ref{H:WAC} holds. Then $\tZ(\la)$ possesses the 
 same property for all $\la\in\Cbb$.
\end{lemma}
\begin{proof}
 Let the Weak Atkinson condition be satisfied for some $\la\in\Cbb$, i.e., whenever $z\equiv \tZ(\la)\.\xi$ 
 on $\Izp$ for some $\xi\in\Cbb^{2n}$ is such that
  \begin{equation}\label{E:definite}
   \sum_0^N z_k^*\.\Ps_k\.z_k=0,
  \end{equation}
 then $z_k\equiv0$ or equivalently $\xi=0$. We claim that the same is true for any $\la\in\Cbb$. By contradiction, 
 suppose that there exist $\nu\in\Cbb\stm\{\la\}$ and $\xi\in\Cbb^{2n}\stm\{0\}$ such that the solution 
 $z_k\coloneq \tZ_k(\nu)\.\xi\not\equiv0$ satisfies equality~\eqref{E:definite}. Then $\Ps_k\.z_k=0$ for all 
 $k\in\Iz$, which together with the equivalent form of system~\eqref{E:Sla} given in~\eqref{E:equiv.Sla} 
 yields that $z$ solves also system~\eqref{E:Sla} on $\Iz$, i.e., $z_k=\Ph_k(\la)\.\zeta$ for all $k\in\Izp$ and 
 some $\zeta\in\Cbb^{2n}\stm\{0\}$, where $\Ph(\la)\coloneq\big(\hZ(\la),\,\,\tZ(\la)\big)$ is the fundamental 
 matrix of system~\eqref{E:Sla} determined by the initial conditions from \eqref{E:Ph.via.al}. This fact implies that 
  \begin{equation*}
   0=\al\.\tZ_0(\nu)\.\xi=\al\.z_0=\al\.\Ph_0(\la)\.\zeta=\big(\al\.\al^*\,\, -\al\.\Jc\.\al^*)\.\zeta=\zeta^{\ssc{[1]}}
  \end{equation*}
 for $\zeta=\big(\zeta^{\ssc{[1]}*}\,\,\,\zeta^{\ssc{[2]}*}\big)^*$ with 
 $\zeta^{\ssc{[1]}},\zeta^{\ssc{[2]}}\in\Cbb^n$, i.e., $z_k=\tZ_k(\la)\.\zeta^{\ssc{[2]}}$ for all $k\in\Izp$. However 
 in that case equality~\eqref{E:definite} contradicts the Weak Atkinson condition for system~\eqref{E:Sla} and the 
 proof is complete.
\end{proof}

If we fix $\la\in\Cbb$, then this hypothesis guarantees that (nontrivial linearly independent combinations of) the 
columns of $\tZ(\la)$ belong to distinct equivalence classes in $\tltp$, i.e., there is at least $n$ linearly 
independent equivalence classes, each of which contains a representative satisfying~\eqref{E:Sla}. This is the main 
difference between the two various Atkinson conditions, because the stronger form ensures the existence of $2n$ 
linearly independent classes for solutions of~\eqref{E:Sla}, see also \cite[Theorem~5.2]{slC.pZ15}. The sufficiency of 
Hypothesis~\ref{H:WAC} stems from the basic characterization of eigenvalues and eigenfunction of system~\eqref{E:Sla}. 
More precisely, it can be easily verified that the eigenvalue problem formulated in~\eqref{E:EP} possesses a nontrivial 
solution if and only if $\la\in\Cbb$ is such that 
 \begin{equation}\label{E:eigenvalue.det.2n}
  \det\mmatrix{-\tZ_{N+1}(\la) & \Jc\.\be^*}=0.
 \end{equation}
This is the discrete counterpart of the classical result for system~\eqref{E:cLHS} given 
in~\cite[Equality~(9.2.9)]{fvA64}, because the boundary conditions from~\eqref{E:EP} are equivalent to the existence of 
a~unique vector $\xi\in\Cbb^{2n}$, for which we have
 \begin{equation}\label{E:boundary.conditions.equiv}
  z_0(\la)=Q^{[0]}\.\xi \qtextq{and} z_{N+1}(\la)=Q^{[N+1]}\.\xi,
 \end{equation}
where $Q^{[0]}\coloneq \big(\Jc\al^*\,\,\,\,0\big)$ and $Q^{[N+1]}\coloneq \big(0\,\,\,\,\Jc\.\be^*\big)$ are 
$2n\times 2n$ matrices such that $Q^{[0]*}\.\Jc\.Q^{[0]}=0=Q^{[N+1]*}\.\Jc\.Q^{[N+1]}$ and 
$\rank\big(Q^{[0]*}\,\,\,\,Q^{[N+1]*}\big)=2n$. Therefore, the $2n\times2n$ matrix appearing
in~\eqref{E:eigenvalue.det.2n} can be written as
 \begin{equation*}
  \mmatrix{-\tZ_{N+1}(\la) & \Jc\.\be^*}=Q^{[N+1]}-\Ph_{N+1}(\la)\.\Ph_0^*(\la)\.Q^{[0]}.
 \end{equation*}
However, the dimension of the latter matrix is unnecessarily large and it can be reduced by half, because 
the eigenvalue problem has, in fact, a nontrivial solution if and only if $\det\be\.\tZ_{N+1}(\la)=0$, see 
\cite[Theorem~2.8]{rSH.pZ14:JDEA}. These nontrivial solutions (i.e., eigenfunctions) are linear combinations of columns 
of $\tZ(\la)$ determined by nonzero vectors $\eta\in\ker\be\.\tZ_{N+1}(\la)$, i.e., the corresponding 
{\it eigenspaces} are
 \begin{equation*}
  \Eigen(\la)\coloneq\big\{\tZ(\la)\.\eta\in\Cbb(\Izp)^{2n} \mid \eta\in\ker\be\.\tZ_{N+1}(\la) \big\},
 \end{equation*}
which means that the value of 
 \begin{equation*}
  0\leq\geom(\la)\coloneq\dim\Eigen(\la)=\dim\ker\be\.\tZ_{N+1}(\la)\leq n
 \end{equation*}
yields the {\it geometric multiplicity} of the eigenvalue $\la$ and it is the same as 
$\dim\ker\big(-\tZ_{N+1}(\la)\,\,\,\,\Jc\.\be^*\big)$. The discreteness of the interval $\Iz$ and the independence 
of the initial value $\tZ_0(\la)$ on $\la$ imply that the function $\det\be\.\tZ_{N+1}(\la)$ is a polynomial in $\la$ 
of degree at most $n(N+1)$. Hence, for any choice of $\al,\be\in\Ga$, system~\eqref{E:Sla} has either only a finite 
number of eigenvalues or the set of all eigenvalues coincides with the whole $\Cbb$, which occurs if and only if 
$\det\be\.\tZ_{N+1}(\la)\equiv0$, see also Example~\ref{Ex:eigenvalue} below. In the latter, very particular, case the 
linear subspace $\Eigen(\la)$ is nontrivial for all $\la\in\Cbb$, while otherwise it holds $\Eigen(\la)\neq\{0\}$ 
only for finitely many (or none) $\la\in\Cbb$ and we have the estimate $\sum_{\la\in\Cbb}\geom(\la)\leq n^2(N+1)$. In 
any case, the first part of Theorem~\ref{T:hilbert} implies the upper bound for the number of linearly independent 
eigenfunctions
 \begin{equation*}
  \dim \bigcup_{\la\in\Cbb}\Eigen(\la)\leq 2n(N+2).
 \end{equation*}
It means that, in any case based only on Hypothesis~\ref{H:basic}, some eigenvalues may be non-real numbers and 
eigenfunctions for distinct eigenvalues may be linearly dependent. Actually, from the 
equivalent form of system~\eqref{E:Sla} given in~\eqref{E:equiv.Sla} one can easily 
deduce that the linear dependence of eigenfunctions $z(\la_1),\dots,z(\la_m)$ corresponding to distinct eigenvalues 
$\la_1,\dots,\la_m\in\Cbb$, i.e., $c_1\.z_k(\la_1)+\dots+c_m\.z_k(\la_m)\equiv0$ on $\Izp$ for some 
$c_1,\dots,c_m\in\Cbb$, is possible only when $c_1\.\la_1\.z_k(\la_1)+\dots+c_m\.\la_m\.z_k(\la_m)\in\ker\Ps_k$ for all 
$k\in\Iz$.

To make the situation more ``self-adjoint'' we need to employ Hypothesis~\ref{H:WAC}. It guarantees that all 
eigenvalues 
are only real numbers, which excludes the possibility $\det\be\.\tZ_{N+1}(\la)\equiv0$, i.e., there is at most $n(N+1)$ 
eigenvalues counting multiplicities. Furthermore, upon applying the Lagrange identity from~\eqref{E:lagrange} with 
$f=g\equiv0$ we obtain that eigenfunctions corresponding to distinct eigenvalues $\la,\nu\in\Rbb$ are linearly 
independent, orthogonal with respect to $\innerP{\cdot}{\cdot}$, and belong to distinct equivalence classes, because 
otherwise
 \begin{equation*}
  0=\innerP{z(\la)}{z(\nu)}=\sum_{k\in\Iz} z_k^*(\la)\.\Ps_k\.z_k(\nu)=\sum_{k\in\Iz} z_k^*(\la)\.\Ps_k\.z_k(\la),
 \end{equation*}
which contradicts the Weak Atkinson condition. Thus,
 \begin{equation}\label{E:sum.geom.estimate1}
  \sum_{\la\in\Cbb}\geom(\la)\leq \sum_{k\in\Iz} \rank \Ps_k\leq (2n-1)\.(N+1).
 \end{equation}
It remains to shed light on the relation between the geometric multiplicity $\geom(\la_j)$ and algebraic multiplicity 
$\alg(\la_j)$ of an eigenvalue $\la_j$, i.e., the multiplicity of $\la_j$ as the root of the polynomial 
$\det\be\.\tZ_{N+1}(\la)$, compare with \cite[Definition~2]{oD.wK07:JDEA}. The polynomial matrix $\be\.\tZ_{N+1}(\la)$ 
can be expressed in the Smith form as the product
 \begin{equation}\label{E:smith}
  \be\.\tZ_{N+1}(\la)=U^{[1]}(\la)\.\diag\{p_1(\la),\dots,p_n(\la)\}\.U^{[2]}(\la)
 \end{equation}
for unique $n\times n$ unimodular polynomial matrices $U^{[1]}(\la),U^{[2]}(\la)$ with 
$\det U^{[1]}(\la)\equiv1\equiv\det U^{[2]}(\la)$ and unique polynomials $p_1(\la),\dots,p_n(\la)\not\equiv0$ such that 
$p_i(\la)$ divides $p_{i+1}(\la)$ for all $i=1,\dots,n-1$, see e.g., \cite[Sections~4.2--4.3]{dsB09}. From the latter 
property it follows that $\geom(\la_j)$ is the same as the number of polynomials $p_i(\la)$ being zero at 
$\la=\la_j$. Therefore, $\geom(\la_j)\leq\alg(\la_j)$ and the estimate given in~\eqref{E:sum.geom.estimate1} can be 
made yet more precise
 \begin{equation*}
  \sum_{\la\in\Cbb} \geom(\la)\leq\min\Big\{n(N+1),\sum_{k\in\Iz}\rank\Ps_k\Big\}.
 \end{equation*}
 
Finally, to establish the equality $\geom(\la_j)=\alg(\la_j)$ for an arbitrary eigenvalue $\la_j$ of~\eqref{E:Sla} we 
need more results from the Weyl--Titchmarsh theory for discrete symplectic system. In \cite[Section~9.5]{fvA64}, 
Atkinson identified a certain {\it characteristic function} as the crucial object in the further study of 
eigenfunctions for linear Hamiltonian differential systems. It is a $2n\times2n$ matrix-valued function appearing in a 
kernel of a unique solution of a nonhomogeneous boundary value problem analogous to~\eqref{E:EE4} below. However, the 
dimension of this function is not in accordance with the case of the second order Sturm--Liouville differential 
equation, where it was convenient to consider a scalar function. We shall rather proceed in a way analogous to the 
latter case, and so our present treatise is based on the $n\times n$ matrix-valued 
{\it Weyl--Titchmarsh $M(\la)$-function}
 \begin{equation*}
  M_{N+1}(\la)=M_{N+1}(\la,\al,\be)\coloneq -[\be\.\tZ_{N+1}(\la)]^{-1}\.\be\.\hZ_{N+1}(\la).
 \end{equation*}
This function is well-defined for all $\la\in\Cbb$ not being an eigenvalue, especially for all $\la\in\Cbb\stm\Rbb$ 
whenever Hypothesis~\ref{H:WAC} holds, and it possesses the reflection property
 \begin{equation}\label{E:EE9}
  M_{N+1}(\la)=M_{N+1}^*(\bla).
 \end{equation}
It gives rise to the so-called {\it Weyl solution} of system~\eqref{E:Sla} given by
 \begin{equation*}
  \Xc(\la)=\Xc(\la,M_{N+1}(\la))\coloneq \Ph(\la)\.\mmatrix{I\\ M_{N+1}(\la)}
            =\hZ(\la)+\tZ(\la)\.M_{N+1}(\la)\in\Cbb(\Izp)^{2n\times n},
 \end{equation*}
which satisfies the boundary conditions $\al\.\Xc_0(\la)=I$ and $\be\.\Xc_{N+1}(\la)=0$. Let us also point out that 
Weyl solutions and $M(\la)$-functions for any $\la,\nu\in\Cbb$ not being eigenvalues with $\nu\neq\bla$ can be 
connected via the semi-inner product as
 \begin{equation}\label{E:Mla.Xla.nu}
  (\bla-\nu)\.\sum_{k\in\Iz} \Xc_k^*(\la)\.\Ps_k\.\Xc_k(\nu)
   \overset{\eqref{E:lagrange}}{=}\Xc_0^*(\la)\.\Jc\.\Xc_0(\nu)=M_{N+1}^*(\la)-M_{N+1}(\nu),
 \end{equation}
because the equality $\be\.\Xc_{N+1}(\la)=0$ implies that all columns of the Weyl solution $\Xc_{N+1}(\la)$ belong 
to $\ker\be=\ran\Jc\be^*$, i.e., $\Xc_{N+1}(\la)=\Jc\be^*\.\Xi(\la)$ for a suitable $\Xi(\la)\in\Cbb^{n\times n}$, 
which yields $\Xc_{N+1}^*(\la)\.\Jc\.\Xc_{N+1}(\nu)=\Xi^*(\la)\.\be\.\Jc\be^*\.\Xi(\nu)=0$. 

In particular, if we take $\nu=\la$, then~\eqref{E:Mla.Xla.nu} reduces to
 \begin{equation}\label{E:EE8}
  0\leq\abs{\im\la}\.\sum_{k\in\Iz} \Xc_k^*(\la)\.\Ps_k\.\Xc_k(\la)=\sgn(\im\la)\.\im M_{N+1}(\la),
 \end{equation}
which together with the holomorphic property of $M_{N+1}(\cdot)$ on $\Cbb\stm\Rbb$ and the reflection property implies, 
under Hypothesis~\ref{H:WAC}, that $M_{N+1}(\cdot):\Cbb\to\Cbb^{n\times n}$ is a {\it Nevanlinna} (or {\it Herglotz} or 
{\it Pick}) {\it function}, see \cite[Section~5]{fG.erT00}. Therefore, it has only isolated singularities, which are 
simple poles located on the real line and coincide with the zeros of $\det\be\.\tZ_{N+1}(\la)$, i.e., the eigenvalues 
of~\eqref{E:Sla}. Hence for every $\la$ in a~suitable neighborhood $\Oc_L(\la_j)$ we have the Laurent series expansion
 \begin{equation}\label{E:EE10}
  M_{N+1}(\la)=(\la-\la_j)^{-1}\.L_{-1}^{[j]}+L_{0}^{[j]}+(\la-\la_j)\.L_{1}^{[j]}+(\la-\la_j)^{2}\.L_{2}^{[j]}+\cdots
 \end{equation}
with coefficients $L_{-1}^{[j]},L_{0}^{[j]},\dots\in\Cbb^{n\times n}$ being Hermitian matrices by~\eqref{E:EE9} and 
$L_{-1}^{[j]}\neq0$. Simultaneously, from~\eqref{E:smith} we can see that the equality $\geom(\la_j)=\alg(\la_j)$ is 
equivalent to the fact that each of polynomials $p_1(\la),\dots,p_n(\la)$ is either constant or has only simple roots, 
i.e.,
 \begin{equation*}
  p_i(\la)=(\la-\la_1)^{m_1}\cdots (\la-\la_r)^{m_r},\quad i=1,\dots,n,
 \end{equation*}
with $\la_1,\dots,\la_r$ being all eigenvalues of~\eqref{E:Sla} and $m_1,\dots,m_r\in\{0,1\}$ being nondecreasing 
as $i$ goes from $1$ to $n$. In other words, it suffices to show that the function 
 \begin{equation*}
  \diag\{1/p_1(\la),\dots,1/p_n(\la)\}=U^{[1]}(\la)\.\big[\be\.\tZ_{N+1}(\la)\big]^{-1} U^{[2]}(\la)
 \end{equation*}
has at $\la=\la_j$, $j=1,\dots,r$, a simple pole. The characterization of eigenvalues of~\eqref{E:Sla} discussed 
above implies that $\be\.\hZ_{N+1}(\la)$ is holomorphic and nonsingular at $\la=\la_j$. So, its inverse 
$\big[\be\.\hZ_{N+1}(\la)\big]^{-1}$ is holomorphic in a suitable neighborhood $\Oc_T(\la_j)$, in which it has a 
certain Taylor series expansion.
Since similar Taylor series expansions can be obtained also for the polynomial matrices
$U^{[1]}(\la)$ and $U^{[2]}(\la)$, it follows that for all $\la\in\Oc_T(\la_j)\cap\Oc_L(\la_j)$ we have
 \begin{align*}
  \diag\{1/p_1(\la),\dots,1/p_n(\la)\}
   &=-U^{[1]}(\la)\.M_{N+1}(\la)\.\big[\be\.\hZ_{N+1}(\la)\big]^{-1}\.U^{[2]}(\la)\\
   &=-U^{[1]}(\la_j)\.L^{[j]}_{-1}\.\big[\be\.\hZ_{N+1}(\la_j)\big]^{-1}\.U^{[2]}(\la_j)\.(\la-\la_j)^{-1}+\dots
 \end{align*}
proving the desired equality, because 
 \begin{equation*}
  \rank U^{[1]}(\la_j)\.L^{[j]}_{-1}\.[\be\.\hZ_{N+1}(\la_j)]^{-1}\.U^{[2]}(\la_j)=\rank L^{[j]}_{-1}\neq0
 \end{equation*} 
and the omitted terms on the right-hand side are multiples of $(\la-\la_j)^r$ for $r\in\{0,1,\dots\}$.

All of these properties are summarized in the following statement, which was partially included in 
Theorem~\ref{T:main.intro}. Let us note that the equality between the algebraic and geometric multiplicities can be 
seen as a natural consequence of the self-adjointness of an underlying linear relation, see 
\cite[Remark~4.4(ii)]{pZ.slC:SAE2} for more details.

\begin{theorem}\label{T:eigenvalues}
 Let matrices $\al,\be\in\Ga$ be given. Then a number $\la\in\Cbb$ is an eigenvalue of~\eqref{E:Sla} if and only if 
 $\det\big(-\tZ_{N+1}(\la)\,\,\,\,\Jc\.\be^*\big)=0$ or equivalently $\det\be\.\tZ_{N+1}(\la)=0$. This is true either 
 for all $\la\in\Cbb$ or only for at most $n(N+1)$ complex numbers $\la$. The geometric multiplicity of every 
 eigenvalue is at most equal to $n$ and all eigenfunctions form a subspace of $\ltp$ such that 
  \begin{equation*}
   \dim\bigcup_{\la\in\Cbb}\Eigen(\la)\leq 2n(N+2).
  \end{equation*}
 If, in addition, Hypothesis~\ref{H:WAC} is satisfied, then there is at most $n(N+1)$ eigenvalues (counting 
 multiplicities), all of which are only real numbers and for each of them the algebraic and geometric multiplicities 
 coincide. Consequently, the eigenvalues can be arranged as
  \begin{equation*}
    \la_1\leq \la_2\leq\dots\leq \la_r<\infty \qtext{with $r\coloneq\sum_{\la\in\Cbb}\alg(\la)\leq n(N+1)$.}
  \end{equation*} 
 Furthermore, eigenfunctions corresponding to distinct eigenvalues are linearly independent, orthogonal with 
 respect to $\innerP{\cdot}{\cdot}$, belong to distinct equivalence classes, and 
   \begin{equation*}
    \dim\bigcup_{\la\in\Cbb}\Eigen(\la)=\sum_{\la\in\Cbb}\geom(\la)=\sum_{\la\in\Cbb}\alg(\la)
                                       \leq\min\Big\{n(N+1),\sum_{k\in\Iz}\rank\Ps_k\Big\}.
  \end{equation*}
\end{theorem}

In the following example we illustrate some of the properties derived above. Let us note that, for simplicity, all 
presented examples concern with systems in the canonical form, i.e., when $\Sc_k\equiv I$ on $\Iz$. This is done 
without 
loss of generality as any system~\eqref{E:Sla} can be transformed in this form by using a suitable unitary 
transformation, see \cite[Remark~2.3]{pZ.slC16}.

\begin{example}\label{Ex:eigenvalue}
 \begin{enumerate}[leftmargin=10mm,topsep=0mm,label={{\normalfont{(\roman*)}}}]
  \item Let $\Iz=\onZ$ be an arbitrary finite discrete interval, a non-decreasing real-valued scalar sequence 
        $v\in\Cbb(\Izp)^1$ be given with $v_0\coloneq0$, and consider system~\eqref{E:Sla} with the coefficients
         \begin{equation*}
          \Sc_k\equiv I_2\in\Cbb^{2\times2} \qtextq{and} \Ps_k=\diag\{0,\De v_k\},
         \end{equation*}
        which obviously satisfy Hypothesis~\ref{H:basic} and yield 
         \begin{equation*}
          \dim\tltp=\sum_{k\in\Iz} \rank \Ps_k=\sum_{k\in\Iz}\sgn(v_{k+1})\leq N+1.
         \end{equation*}
        If we take $\al=(0\,\,\,\,1)$, then the corresponding fundamental matrix is
         \begin{equation*}
          \Ph_k(\la)=\mmatrix{\la\.v_k & -1\\ 1 & 0} \qtext{for all $\la\in\Cbb$,}
         \end{equation*}
        and so the Weak Atkinson condition is never true. The choice $\be=(0\,\,\,\,1)$ leads to
         \begin{equation*}
          \det\be\.\tZ_{N+1}(\la)=\mmatrix{0 & 1}\mmatrix{-1\\ 0}=0, 
         \end{equation*}
        which means that every $\la\in\Cbb$ is a simple eigenvalue with the unique (up to a nonzero constant multiple) 
        eigenfunction $z_k(\la)\equiv\msmatrix{-1\\ 0}$ belonging to the zero equivalence class $[0]$. On the other 
        hand, the choice $\be=(1\,\,\,\,0)$ leads to
         \begin{equation*}
          \det\be\.\tZ_{N+1}(\la)=\mmatrix{1 & 0}\mmatrix{-1\\ 0}=-1, 
         \end{equation*}
        which means that there is no eigenvalue.
        
        Similarly, if we take $\al=(1\,\,\,\,0)$, then we have the fundamental matrix
         \begin{equation*}
          \Ph_k(\la)=\mmatrix{1 & \la\.v_k\\ 0 & 1}
         \end{equation*}
        and the Weak Atkinson condition is satisfied whenever $v_k>0$ for some $k\in\Izp\stm\{0\}$. In that case, for 
        $\be=(0\,\,\,\,1)$ there is no eigenvalue, while for $\be=(1\,\,\,\,0)$ we have the only eigenvalue $\la=0$ 
        with the unique (up to a nonzero constant multiple) eigenfunction $z_k(0)\equiv\msmatrix{0\\ 1}$ with 
        the norm $\normP{z(0)}=\sqrt{v_{N+1}}\neq0$, i.e., $z(0)\not\in[0]$.
        
  \item Let $\Iz=\onZ$ be an arbitrary finite discrete interval and consider system~\eqref{E:Sla} with
         \begin{equation*}
          \Sc_k\equiv I_4\in\Cbb^{4\times4} \qtextq{and} 
          \Ps_k\equiv\mmatrix{a\.I_2 & \sqrt{ab}\.I_2\\ \sqrt{ab}\.I_2 & b\.I_2},
         \end{equation*}
        where $a,b>0$ are given numbers. Then Hypothesis~\ref{H:basic} is satisfied with 
         \begin{equation*}
          \Vc_k\equiv\msmatrix{-\sqrt{ab}\.I_2 & -b\.I_2\\ a\.I_2 & \sqrt{ab}\.I_2} \qtextq{and}
          \dim\tltp=\sum_{k\in\Iz} \rank \Ps_k=2(N+1). 
         \end{equation*}
        If we choose $\al=(I_2\,\,\,\,0)$, then we obtain the fundamental matrix
         \begin{equation*}
          \Ph_k(\la)=(I+\la\.\Vc_k)^{-k}=\mmatrix{I_2+k\.\la\.\sqrt{ab}\.I_2 & k\.\la\.b\.I_2\\
                                                  -k\.\la\.a\.I_2            & I_2-k\.\la\.\sqrt{ab}\.I_2}.
         \end{equation*}
        Since in that case $\tZ_k^*(\la)\.\Ps_k\tZ_k(\la)\equiv b\.I_2$, the positivity of $b$ guarantees that 
        Hypothesis~\ref{H:WAC} holds. For $\be=(I_2\,\,\,\,0)$ we have the double eigenvalue $\la=0$ with the constant 
        pair of linearly independent eigenfunctions
         \begin{equation*}
          z^{[1]}_k\equiv \mmatrix{0 & 0 & 1 & 0}^{\!\top} \qtextq{and} 
          z^{[2]}_k\equiv \mmatrix{0 & 0 & 0 & 1}^{\!\top}.
         \end{equation*}
        Similarly, for $\be=(0\,\,\,\,I_2)$ we have the double eigenvalue $\la=1/[(N+1)\.\sqrt{ab}]$ with the pair of 
        linearly independent eigenfunctions
         \begin{equation*}
          z^{[1]}_k=\mmatrix{\frac{kb}{(N+1)\.\sqrt{ab}} & 0 & 1-\frac{k}{N+1} & 0}^{\!\!\top} \qtextq{and}
          z^{[2]}_k=\mmatrix{0 & \frac{kb}{(N+1)\.\sqrt{ab}} & 0 & 1-\frac{k}{N+1}}^{\!\!\top}\!\!  
         \end{equation*} 
        for all $k\in\Izp$. Finally, for $\be=\msmatrix{1 & 0 & 0 & 0\\ 0 & 0 & 0 & 1}$ we have the simple eigenvalues
        $\la_1=0$ and $\la_2=1/[(N+1)\.\sqrt{ab}]$ with the corresponding eigenfunctions
         \begin{equation*}
          z_k(\la_1)\equiv \mmatrix{0 & 0 & 1 & 0}^{\!\top} \qtextq{and}
          z_k(\la_2)=\mmatrix{0 & \frac{kb}{(N+1)\.\sqrt{ab}} & 0 & 1-\frac{k}{N+1}}^{\!\!\top}, \qtext{$k\in\Izp$.} 
         \end{equation*}  
        All of these eigenfunctions satisfy
         \begin{equation*}
          \normP{z^{[1]}}=\normP{z^{[2]}}=\normP{z(\la_1)}=\normP{z(\la_2)}=\sqrt{b(N+1)} 
         \end{equation*}
        and one can easily verify the mutual orthogonality of the eigenfunctions in the pairs as well as the fact that 
        each eigenfunction in the pair belongs to a distinct equivalence class in $\tltp$.
 \end{enumerate}
\end{example}

While eigenfunctions corresponding to distinct eigenvalues can be easily orthonormalized due to 
Hypothesis~\ref{H:WAC}, linearly independent eigenfunctions corresponding to the same eigenvalue can be orthonormalized 
by using the Gram--Schmidt process as follows. Let $\la_j\in\Rbb$ be an eigenvalue of system~\eqref{E:Sla} with 
$\geom(\la_j)$ linearly independent eigenfunctions, i.e., it holds
 \begin{equation*}
  1\leq r_j\coloneq\geom(\la_j)=\dim\ker\be\.\tZ_{N+1}(\la_j)\leq n
 \end{equation*}
and there are linearly independent vectors
 \begin{equation}\label{E:EE0}
  \eta^{[1]},\dots,\eta^{[r_j]}\in\ker\be\.\tZ_{N+1}(\la_j)\subseteq\Cbb^{n},
 \end{equation}
which yield all linearly independent eigenfunctions
 \begin{equation*}
  z^{[\ell]}=z^{[\ell]}(\la_j)\coloneq \tZ(\la_j)\.\eta^{[\ell]} \qtext{for all $\ell\in\{1,\dots,r_j\}$.}
 \end{equation*}
However these eigenfunctions may not be orthogonal, so we aim now to find a basis of $\ker\be\.\tZ_{N+1}(\la_j)$ with 
the property
 \begin{equation*}
  \innerP{z^{[j]}}{z^{[\ell]}}=\eta^{[j]*}\.\sum_{s=0}^N \tZ_s^*(\la_j)\.\Ps_s\.\tZ_s(\la_j)\.\eta^{[\ell]}
   =\begin{cases}
     0, & j\neq\ell,\\[1mm]
     1, & j=\ell.
    \end{cases}
 \end{equation*}
Let us start with an arbitrary orthonormal basis of $\ker\be\.\tZ_{N+1}(\la_j)$ with respect to the standard inner 
product in $\Cbb^n$, i.e., we choose $\xi^{[1]},\dots,\xi^{[r_j]}\in\ker\be\.\tZ_{N+1}(\la_j)$ so that
 \begin{equation*}
  \xi^{[j]*}\.\xi^{[\ell]}   
   =\begin{cases}
     0, & j\neq\ell,\\[1mm]
     1, & j=\ell.
    \end{cases}
 \end{equation*}
If we denote the $n\times n$ matrix 
 \begin{equation}\label{E:Om.def}
  \Om\coloneq \sum_{k\in\Iz} \tZ_k^*(\la_j)\.\Ps_k\.\tZ_k(\la_j),
 \end{equation}
then our original goal reads as to find vectors from~\eqref{E:EE0} such that
 \begin{equation}\label{E:EE00}
  \eta^{[j]*}\.\Om\,\eta^{[\ell]}   
   =\begin{cases}
     0, & j\neq\ell,\\[1mm]
     1, & j=\ell.
    \end{cases}
 \end{equation}
Since $\Om>0$ by Hypothesis~\ref{H:WAC}, the square root $\Om^{1/2}$ and its inverse are well-defined, so the 
vectors
 \begin{equation*}
  \om^{[1]}\coloneq \Om^{1/2}\.\xi^{[1]},\dots,\om^{[r_j]}\coloneq \Om^{1/2}\.\xi^{[r_j]}
 \end{equation*}
are linearly independent and they form a basis of an $r_j$-dimensional subspace $\Cbb^n$. This subspace possesses 
also an orthonormal basis $\rh^{[1]},\dots,\rh^{[r_j]}\in\Cbb^n$, which gives rise to the vectors
 \begin{equation}\label{E:EE0b}
  \eta^{[1]}\coloneq \Om^{-1/2}\.\rh^{[1]},\dots,\eta^{[r_j]}\coloneq \Om^{-1/2}\.\rh^{[r_j]}.
 \end{equation}
Then for all $j\in\{1,\dots,r_j\}$ and suitable numbers $c_1,\dots,c_{r_j}\in\Cbb$ we can write
 \begin{align*}
  \eta^{[j]}=\Om^{-1/2}\.\rh^{[j]}
            &=\Om^{-1/2}\.(c_1\.\om^{[1]}+\dots+c_{r_j}\.\om^{[r_j]})\\
            &=c_1\.\xi^{[1]}+\dots+c_{r_j}\.\xi^{[r_j]}\in\ker\be\.\tZ_{N+1}(\la_j)
 \end{align*}
and, at the same time, these vectors obviously satisfy~\eqref{E:EE00}, which means that $\eta^{[1]},\dots,\eta^{[r_j]}$ 
yield an orthonormal set of eigenfunctions corresponding to the eigenvalue $\la_j$. This construction justifies the 
following statement.

\begin{theorem}\label{T:orthonormal.set}
 Let $\al,\be\in\Ga$ be given and Hypothesis~\ref{H:WAC} hold. Then system~\eqref{E:Sla} possesses a finite orthonormal 
 set $z^{[1]}(\la_1),\dots,z^{[r_1]}(\la_1),z^{[1]}(\la_2),\dots\in\Cbb(\Izp)^{2n}$ of all eigenfunctions.
\end{theorem}

This orthonormal set shall play a crucial role in our main goal concerning a generalization of the result concerning 
the eigenfunction expansion on an arbitrary finite discrete interval $\Iz$ recalled in Theorem~\ref{T:EE.admissible}. 
It is also closely connected with the residue of the $M(\la)$-function. To show that, we need to derive an 
explicit formula for the solution of the nonhomogeneous problem
 \begin{equation}\label{E:EE4}
  \Slaf{\la}{f},\quad \la\in\Cbb,\quad k\in\Iz,\quad f\in\Cbb(\Izp),\quad \al\.z_0=0=\be\.z_{N+1}.
 \end{equation}
Let us start with the classical variation of constants formula, i.e., assume that the solution of~\eqref{E:EE4} is of 
the form $\hz(\la)=\Ph(\la)\.\xi$ for a suitable $\xi\in\Cbb(\Izp)^{2n}$. Upon inserting $\hz$ into~\Slaf{\la}{f} we get
 \begin{equation*}
  \hz_k(\la)=\Ph_k(\la)\.\xi_{k+1}-\Jc\.\Ps_k\.f_k
 \end{equation*}
for all $k\in\Iz$, which immediately yields
 \begin{equation*}
  \De\xi_k=\Ph^{-1}_k(\la)\.\Jc\.\Ps_k\.f_k \qtextq{or equivalently}
  \xi_k=\xi_0+\sum_{j=0}^{k-1} \Ph^{-1}_j(\la)\.\Jc\.\Ps_j\.f_j.
 \end{equation*}
Since the boundary condition at $k=0$ implies
 \begin{equation*}
  0=\al\.\Ph_0(\la)\.\xi_0=(I,0)\.\xi_0,\qtextq{i.e.,} \xi_0=\mmatrix{0\\ \xi^{[2]}_0},
 \end{equation*}
we obtain the expression
 \begin{equation*}
  \hz_k(\la)=\tZ_k(\la)\.\xi^{[2]}_0+\sum_{j=0}^{k-1} \Ph_k(\la)\.\Jc\Ph_j^*(\bla)\.\Ps_j\.f_j.
 \end{equation*}
If $\la$ is not an eigenvalue, then from the boundary condition $k=N+1$ and the invertibility of the matrix
$\be\.\tZ_{N+1}(\la)$ it follows
 \begin{equation*}
  \xi^{[2]}_0=\big(M_{N+1}(\la),-I\big)\.\sum_{j\in\Iz} \Jc\.\Ph_j^*(\bla)\.\Ps_j\.f_j,
 \end{equation*}
and so we have 
 \begin{equation}\label{E:EE5}
  \hz_k(\la)=\tZ_k(\la)\.\sum_{j\in\Iz}\big[M_{N+1}(\la)\.\tZ_j^*(\bla)+\hZ_j^*(\bla)\big]\.\Ps_j\.f_j
             +\sum_{j=0}^{k-1} \Ph_k(\la)\.\Jc\.\Ph_j^*(\bla)\.\Ps_j\.f_j.
 \end{equation}
Therefore, we get the following result, which does not rely on Hypothesis~\ref{H:WAC}, cf. \cite[Section~4]{slC.pZ15}.

\begin{theorem}\label{T:T5+C5A}
 Let $\al,\be\in\Ga$ and $f\in\Cbb(\Izp)^{2n}$ be given. Whenever $\la\in\Cbb$ is not an eigenvalue of~\eqref{E:Sla},
 the sequence $\hz\in\Cbb(\Izp)^{2n}$ from~\eqref{E:EE5} represents the unique solution of problem~\eqref{E:EE4}. 
 Moreover, the unique solution of the modified boundary value problem
  \begin{equation*}
   \Slaf{\la}{f},\quad \la\in\Cbb,\quad k\in\Iz,\quad f\in\Cbb(\Izp),\quad \al\.z_0=\xi,\quad \be\.z_{N+1}=0
  \end{equation*}
 for a given $\xi\in\Cbb^n$ can be expressed by using the Weyl solution as
  \begin{equation*}
   z(\la)=\hz(\la)+\Xc(\la)\.\xi \qtext{on $\Izp$.}
  \end{equation*}
\end{theorem}
\begin{proof}
 It follows from the previous calculations. We only note that the uniqueness of the solution of~\eqref{E:EE4} is 
 guaranteed by the fact that $\la$ is not an eigenvalue, because otherwise the difference of two solutions, being 
 nontrivial, would satisfy the eigenvalue problem from~\eqref{E:EP}, i.e., the number $\la$ would be an eigenvalue. For 
 the second part we utilize the facts that $\al\.\Xc_0(\la)=I$ and $\be\.\Xc_{N+1}(\la)=0$.
\end{proof}

Splitting the first sum in~\eqref{E:EE5} in two parts we get
 \begin{equation*}
  \begin{aligned}
   \hz_k=\sum_{j=0}^{k-1} \Big\{\tZ_k(\la)\.\big[M_{N+1}(\la)\.\tZ^*_j(\bla)+\hZ_j^*(\bla)\big]
             &+\Ph_k(\la)\.\Jc\.\Ph_j^*(\bla)\Big\}\.\Ps_j\.f_j\\
             &\hspace{-20mm}+\sum_{j=k}^N \tZ_k(\la)\.\big[M_{N+1}(\la)\.\tZ^*_j(\bla)+\hZ_j^*(\bla)\big]\.\Ps_j\.f_j,
  \end{aligned}
 \end{equation*}
which leads to
 \begin{equation}\label{E:EE5a}
  \hz_k(\la)=\sum_{j\in\Iz} G_{k,j}(\la)\.\Ps_j\.f_j, \quad k\in\Izp,
 \end{equation}
with the {\it kernel} (or {\it Green function}) $G_{k,j}(\la)$ being defined as
 \begin{equation}\label{E:EE5b}
  G_{k,j}(\la)\coloneq \begin{cases}
                        \Xc_k(\la)\.\tZ^*_j(\bla), & j\in[0,k-1]_\sZbb,\\[1mm]
                        \tZ_k(\la)\.\Xc_j^*(\bla), & j\in[k,N+1]_\sZbb,
                       \end{cases}
 \end{equation}
compare with \cite[Formula~(4.5)]{slC.pZ15}. Since the term $\Ph_k(\la)\.\Jc\.\Ph_j^*(\bla)$ is not involved in the 
second sum and $\Ph_k(\la)\.\Jc\.\Ph_k^*(\bla)\equiv\Jc$, we may say, in a certain sense, that the function 
$G_{k,j}(\la)$ has the jump equal to $-\Jc$ as $j$ growths from $0$ to $N$ through the index $k$. It is also quite 
straightforward to verify that 
 \begin{equation*}
  G_{k,j}^*(\la)= \begin{cases}
                   G_{j,k}(\bla), & j\neq k,\\[1mm]
                   G_{k,k}(\bla)+\Jc, & j=k.
                  \end{cases}
 \end{equation*}
Furthermore, if $\nu\in\Cbb$ is not an eigenvalue and $y(\nu)$ denotes the unique solution of yet another 
nonhomogeneous problem
 \begin{equation*}
  \Slaf{\nu}{g},\quad k\in\Iz,\quad g\in\Cbb(\Izp),\quad \al\.z_0=0=\be\.z_{N+1},
 \end{equation*}
then equality~\eqref{E:lagrange} and the equivalent form of the boundary conditions displayed 
in~\eqref{E:boundary.conditions.equiv} yield
 \begin{equation*}
  0=\hz_k^*(\la)\.\Jc\.y_k(\nu)\Big|_{0}^{N+1}
   =\sum_{k\in\Iz}\Big\{(\bla-\nu)\.\hz_k^*(\la)\.\Ps_k\.y_k(\nu)+f_k^*\.\Ps_k\.y_k(\nu)-z_k^*(\la)\.\Ps_k\.g_k\Big\},
 \end{equation*}
which upon substituting the expression of $\hz(\la)$ from~\eqref{E:EE5a} and its analogue for $y(\nu)$ 
leads to 
 \begin{equation*}
  \sum_{j\in\Iz} \sum_{k\in\Iz} f_j^*\.\Ps_j\Big\{G_{k,j}^*(\la)-G_{j,k}(\nu)
                     -(\bla-\nu)\sum_{s\in\Iz} G_{s,j}^*(\la)\.\Ps_s\.G_{s,k}(\nu)\Big\}\.\Ps_k\.g_k=0.
 \end{equation*}
More generally, if $j\in\Iz$ is fixed, then the kernel $G_{k,j}(\la)$ satisfies system~\eqref{E:Sla} for all 
$k\in\Iz\stm\{j\}$, while for $k=j$ we have $G_{k,k}(\la)-\Sbb_k(\la)\.G_{k+1,k}(\la)=-\Jc$. Therefore, if we 
replace the nonhomogeneous term $\Jc\.\Ps\.f$ on the right-hand side of~\eqref{E:Slaf} by $\Jc f$, then 
$G_{\cdot,j}(\la)$ solves the matrix extension of this modified system on $\Iz$ with $f\in\Cbb(\Izp)^{2n\times2n}$ 
being such that $f_j=I$ and $f_k=0$ for $k\in\Izp\stm\{j\}$. In that case, from the analogue of the Extended Lagrange 
identity for this modified system we get the following statement, compare with \cite[Theorems~9.4.2--9.4.3]{fvA64}.

\begin{lemma}
 Let $\al,\be\in\Ga$ be given. If the complex numbers $\la,\nu\in\Cbb$ are not eigenvalues, then the kernel 
 $G_{\cdot,\cdot}(\cdot)$ satisfies
  \begin{align*}
   G_{k,j}^*(\la)-G_{j,k}(\nu)=&(\bla-\nu)\sum_{s\in\Iz} G_{s,j}^*(\la)\.\Ps_s\.G_{s,k}(\nu)
                               +(\bla-\nu)\.G_{k,j}^*(\la)\.\Ps_k\.\Jc\\
                               &\hspace*{7mm}-(\bla-\nu)\.\Jc\.\Ps_j\.G_{j,k}^*(\nu)
                                +\big(1-\abs{\sgn(j-k)}\big)\.\Jc\\
                               &\hspace*{15mm}-(\bla-\nu)\.\big(1-\abs{\sgn(j-k)}\big)\.\Jc\.\Ps_j\.\Jc
  \end{align*}
 for any pair $k,j\in\Iz$. In particular, 
  \begin{equation*}
   \Ps_j\.\big[G_{k,j}^*(\la)-G_{j,k}(\nu)\big]\.\Ps_k
    =(\bla-\nu)\sum_{s\in\Iz} \Ps_j\.G_{s,j}^*(\la)\.\Ps_s\.G_{s,k}(\nu)\.\Ps_k.
  \end{equation*}
\end{lemma}

Now we can calculate the residue $L_{-1}^{[j]}$ from the Laurent expansion of the $M(\la)$-function in a~neighborhood 
of an eigenvalue $\la_j$ used in~\eqref{E:EE10}. In general, it is equal to a Hermitian matrix obtained as the 
(finite) limit
 \begin{equation*}
  L_{-1}^{[j]}=\lim_{\la\to\la_j} (\la-\la_j)\.M_{N+1}(\la)
 \end{equation*}
and, according to the theory of Nevanlinna matrix-valued functions, it is even a~negative semidefinite matrix, see 
\cite[Theorem~5.4(vii)]{fG.erT00}. This fact is justified by the following theorem, in which we derive the 
explicit form of the residue matrix.
 
\begin{theorem}\label{T:T7}
 Let $\al,\be\in\Ga$ be given and Hypothesis~\ref{H:WAC} hold. The residue of $M_{N+1}(\la)$ at any eigenvalue 
 $\la_j$ of~\eqref{E:Sla} is 
  \begin{equation*}
   L_{-1}^{[j]}=-\sum_{\ell=1}^{r_j} \eta_j^{[\ell]}\.\eta_j^{[\ell]*} \qtextq{and}
   \rank L_{-1}^{[j]}\leq r_j\coloneq\dim\ker\be\.\tZ_{N+1}(\la_j)
  \end{equation*}
 for $\eta_j^{[\ell]}\in\ker\be\.\tZ_{N+1}(\la_j)$ being the vectors constructed in~\eqref{E:EE0b}. Consequently, the 
 residue of $G_{k,s}(\la)$ at $\la_j$ is equal to
  \begin{equation*}
   -\sum_{\ell=1}^{r_j} z_k^{[\ell]}(\la_j)\.z_s^{[\ell]*}(\la_j).
  \end{equation*}
\end{theorem}
\begin{proof}
 Let $\la_j$ be an arbitrary eigenvalue of~\eqref{E:Sla}. First of all, we show the inequality 
 $\rank L_{-1}^{[j]}\leq r_j$. By definition, we have
  \begin{equation}\label{E:EE11a}
   \be\.\tZ_{N+1}(\la)\.M_{N+1}(\la)=-\be\.\hZ_{N+1}(\la)
  \end{equation}
 with the right hand-side being holomorphic for all $\la\in\Cbb$. Upon inserting~\eqref{E:EE11a} into the 
 Laurent expansion of $M_{N+1}(\la)$ given in~\eqref{E:EE10} we get
  \begin{equation*}
   (\la-\la_j)^{-1}\.\be\.\tZ_{N+1}(\la)\.L_{-1}^{[j]}+\be\.\tZ_{N+1}(\la)\.L_{0}^{[j]}
     +(\la-\la_j)\.\be\.\tZ_{N+1}(\la)\.L_{1}^{[j]}+\dots
    =-\be\.\hZ_{N+1}(\la)
  \end{equation*}
 in a punctured neighborhood $\Oc^*(\la_j)\coloneq\Oc(\la_j)\stm\{\la_j\}$, so
  \begin{equation*}
   \be\.\tZ_{N+1}(\la)\.L_{-1}^{[j]}(\la-\la_j)^{-1}
  \end{equation*}
 is bounded in $\Oc^*(\la_j)$, i.e., $\abs{(\la-\la_j)^{-1}}\times \norm{\be\.\tZ_{N+1}(\la)\.L_{-1}^{[j]}}\leq K$ for 
 all $\la\in\Oc^*(\la_j)$, some $K>0$, and any matrix norm $\norm{\cdot}$, which is possible only when
  \begin{equation*}
   \be\.\tZ_{N+1}(\la_j)\.L_{-1}^{[j]}=0.
  \end{equation*}
 Thus, by Sylvester's rank inequality, 
 $0=\rank\be\.\tZ_{N+1}(\la)\.L_{-1}^{[j]}\geq \rank\be\.\tZ_{N+1}(\la_j)+\rank L_{-1}^{[j]}-n$, 
 which yields
  \begin{equation}\label{E:residue.rank.inequality}
   \rank L_{-1}^{[j]}\leq r_j.
  \end{equation}
 
 The kernel $G_{k,j}(\la)$ is holomorphic for all $\la\in\Cbb$ except for the eigenvalues of~\eqref{E:Sla} being its 
 poles. Upon inserting the expression~\eqref{E:EE10} into the definition of $G_{k,j}(\la)$ given in~\eqref{E:EE5b} we 
 obtain
  \begin{equation}\label{E:EE13}
   G_{k,j}(\la)=\begin{cases}
                 \tZ_k(\la)\.\big[(\la-\la_j)^{-1}\.L_{-1}^{[j]}+L_{0}^{[j]}+\cdots\big]\.\tZ_j^*(\bla)
                              +\hZ_k(\la)\.\tZ^*_j(\bla),&\!\! j\in[0,k-1]_\sZbb,\\[1mm]
                 \!\tZ_k(\la)\.\big[(\la-\la_j)^{-1}\.L_{-1}^{[j]}+L_{0}^{[j]}+\cdots\big]\.\tZ_j^*(\bla) 
                              +\tZ_k(\la)\.\hZ^*_j(\bla),&\!\! j\in[k,N+1]_\sZbb, 
                \end{cases}
  \end{equation}
 for all $\la\in\Oc^*(\la_j)$. Since every eigenfunction $z^{[\ell]}(\la_j)$, $\ell\in\{1,\dots,r_j\}$, 
 from the orthonormal set established in Theorem~\ref{T:orthonormal.set} can be seen as a nontrivial solution of the 
 nonhomogeneous problem
  \begin{equation*}
   z_k^{[\ell]}(\la_j)=\Sbb(\la)\.z_{k+1}^{[\ell]}(\la_j)-(\la_j-\la)\.\Jc\.\Ps_k\.z_k^{[\ell]}(\la_j),\quad
   k\in\Iz,\quad \al\.z_0^{[\ell]}(\la_j)=0=\be\.z_{N+1}^{[\ell]}(\la_j),
  \end{equation*}
 Theorem~\ref{T:T5+C5A} and equality~\eqref{E:EE5a} imply that
  \begin{equation}\label{E:EE13a}
   z_k^{[\ell]}(\la_j)=(\la_j-\la)\.\sum_{s\in\Iz} G_{k,s}(\la)\.\Ps_s\.z_s^{[\ell]}(\la_j)
  \end{equation}
 for any $\la\in\Cbb$ not being an eigenvalue. Consequently, this equality together~\eqref{E:EE13} yields
  \begin{align*}
   z_k^{[\ell]}(\la_j)
     &=-\sum_{s=0}^{k-1}\Big\{\tZ_k(\la)\.\big[L_{-1}^{[j]}+(\la-\la_j)\.L_{0}^{[j]}+\cdots\big]\.\tZ_s^*(\bla)
              +(\la-\la_j)\.\hZ_k(\la)\.\tZ^*_l(\bla)\Big\}\.\Ps_s\.z_s^{[\ell]}(\la_j)\\
      &\hspace*{6mm}
       -\sum_{s=k}^{N}\Big\{\tZ_k(\la)\.\big[L_{-1}^{[j]}+(\la-\la_j)\.L_{0}^{[j]}+\cdots\big]\.\tZ_s^*(\bla)
              +(\la-\la_j)\.\tZ_k(\la)\.\hZ^*_l(\bla)\Big\}\.\Ps_s\.z_s^{[\ell]}(\la_j),
  \end{align*}
 and so as $\la$ tends to $\la_j$ we get
  \begin{equation*}
   z_k^{[\ell]}(\la_j)=-\tZ_k(\la_j)\.L_{-1}^{[j]}\sum_{s\in\Iz} \tZ_s^*(\la_j)\.\Ps_s\.z_s^{[\ell]}(\la_j)
  \end{equation*}
 or equivalently
  \begin{equation*}
   \eta_j^{[\ell]}=-L_{-1}^{[j]}\sum_{s\in\Iz} \tZ_s^*(\la_j)\.\Ps_s\.\tZ_s(\la_j)\.\eta_j^{[\ell]}
   \qtextq{or} \eta_j^{[\ell]}=-L_{-1}^{[j]}\.\Om\,\eta_j^{[\ell]}
  \end{equation*}
 for the corresponding $\eta_j^{[\ell]}\in\ker\be\.\tZ_{N+1}(\la_j)$ and the positive definite matrix $\Om$ introduced 
 in~\eqref{E:Om.def}. The latter relation can be written also as
  \begin{equation*}
   \Om^{1/2}\.\eta_j^{[\ell]}=-(\Om^{1/2}\.L_{-1}^{[j]}\.\Om^{1/2})\.\Om^{1/2}\.\eta_j^{[\ell]}, \qtextq{i.e.,}
   \rh_j^{[\ell]}=-(\Om^{1/2}\.L_{-1}^{[j]}\.\Om^{1/2})\.\rh_j^{[\ell]}
  \end{equation*}
 by~\eqref{E:EE0b}, which shows that $-\Om^{1/2}\.L_{-1}^{[j]}\.\Om^{1/2}$ acts as the identity operator on the 
 subspace spanned by the orthonormal set of vectors $\rh_j^{[1]},\dots,\rh_j^{[r_j]}$, i.e., 
 $\rank \Om^{1/2}\.L_{-1}^{[j]}\.\Om^{1/2}\geq r_j$. Therefore, the invertibility of $\Om$
 together with~\eqref{E:residue.rank.inequality} yields
  \begin{equation*}
   \rank \Om^{1/2}\.L_{-1}^{[j]}\.\Om^{1/2}=r_j,
  \end{equation*}
 i.e., the matrix $\Om^{1/2}\.L_{-1}^{[j]}\.\Om^{1/2}$ has $r_j$ eigenvalues equal to $1$ and $n-r_j$ eigenvalues equal 
 to $0$. At the same time, $-\Om^{1/2}\.L_{-1}^{[j]}\.\Om^{1/2}$ is an orthogonal projector matrix onto 
 $r_j$-dimensional subspace $\Cbb^n$ spanned by the vectors in~\eqref{E:EE0b}. Thus, 
  \begin{equation*}
   -\Om^{1/2}\.L_{-1}^{[j]}\.\Om^{1/2}=\Om^{1/2}\.\sum_{\ell=1}^{r_j} \eta_j^{[\ell]}\.\eta_j^{[\ell]*}\.\Om^{1/2},
  \end{equation*}
 which proves the first part of the statement. Consequently, by~\eqref{E:EE13} we have
  \begin{equation}\label{E:T7.proof.last.step}
   G_{k,s}(\la)=
   -(\la-\la_j)^{-1}\.\tZ_k(\la)\.\sum_{\ell=1}^{r_j} \eta_j^{[\ell]}\.\eta_j^{[\ell]*}\.\tZ_s^*(\la)+\dots,
  \end{equation}
 with the omitted terms on the right-hand side of~\eqref{E:T7.proof.last.step} being multiples of $(\la-\la_j)^m$ for 
 $m=0,1,\dots$, which yields that
  \begin{equation*}
   \lim_{\la\to\la_j} (\la-\la_j)\.G_{k,s}(\la)=-\sum_{\ell=1}^{r_j} z_k^{[\ell]}(\la_j)\.z_s^{[\ell]*}(\la_j)
  \end{equation*}
 as stated in the second part of the theorem. 
\end{proof}

The last result of this section shows yet another connection between eigenfunctions and the kernel 
$G_{\cdot,\cdot}(\cdot)$, which turns out to be useful in a further study of the $M(\la)$-function.

\begin{lemma}\label{L:L7A}
 Let $\al,\be\in\Ga$ be given and Hypothesis~\ref{H:WAC} hold. Whenever $\la\in\Cbb$ is not an eigenvalue it holds
  \begin{equation}\label{E:EE14b}
   \sum_{j=1}^{r}\sum_{\ell=1}^{r_j} \abs{\la-\la_j}^{-2} z_k^{[\ell]}(\la_j)\.z_k^{[\ell]*}(\la_j)
    \leq \sum_{s\in\Iz} G_{k,s}(\la)\.\Ps_s\.G^*_{k,s}(\la)
  \end{equation}
 for $r$ being as in Theorem~\ref{T:eigenvalues} and $k\in\Izp$. Consequently, we have
  \begin{equation*}
   \sum_{j=1}^{r}\sum_{\ell=1}^{r_j} \abs{\la-\la_j}^{-2} z_k^{[\ell]*}(\la_j)\.z_k^{[\ell]}(\la_j)
    \leq \tr\sum_{s\in\Iz} G_{k,s}(\la)\.\Ps_s\.G^*_{k,s}(\la).
  \end{equation*}
 Especially, for $\la=i$ we obtain
  \begin{equation*}
   \sum_{j=1}^{r}\sum_{\ell=1}^{r_j} \frac{z_k^{[\ell]*}(\la_j)\.z_k^{[\ell]}(\la_j)}{1+\la_j^2}
    \leq \tr\sum_{s\in\Iz} G_{k,s}(i)\.\Ps_s\.G^*_{k,s}(i).
  \end{equation*}
\end{lemma}
\begin{proof}
 The first estimate follows by direct calculation from the inequality
  \begin{equation}\label{E:L7A.full.product}
   \begin{aligned}
    &\sum_{s\in\Iz} \Big[G_{k,s}(\la)
      +\sum_{j=1}^{r}\sum_{\ell=1}^{r_j} (\la-\la_j)^{-1} z_k^{[\ell]}(\la_j)\.z_s^{[\ell]*}(\la_j)\Big]
        \times\Ps_s\times\\
          &\hspace*{25mm}\times\Big[G_{k,s}(\la)
             +\sum_{j=1}^{r}\sum_{\ell=1}^{r_j} (\la-\la_j)^{-1} z_k^{[\ell]}(\la_j)\.z_s^{[\ell]*}(\la_j)\Big]^*\geq0,
   \end{aligned}
  \end{equation}
 the expression of $z^{[\ell]}(\la_j)$ given in~\eqref{E:EE13a}, and the orthonormality of the system of 
 eigenfunctions. More precisely, if we abbreviate the double sum 
 $\sum_{j=1}^{r}\sum_{\ell=1}^{r_j} z_k^{[\ell]}(\la_j)\.z_s^{[\ell]*}(\la_j)$ as 
 $\sum_{\la_j} z_k^{[j]}\.z_s^{[j]*}$ for simplicity, then after expanding the left-hand side 
of~\eqref{E:L7A.full.product} we get
  \begin{align*}
   0&\leq\sum_{s\in\Iz} G_{k,s}(\la)\.\Ps_s\.G_{k,s}^*(\la)
     +\sum_{s\in\Iz}\sum_{\la_j} (\la-\la_j)^{-1}\.z_k^{[j]}\.z_s^{[j]*}\Ps_s\.G_{k,s}^*(\la)\\
    &\hspace*{25mm}+\sum_{s\in\Iz}\sum_{\la_j} (\bla-\la_j)^{-1}\.G_{k,s}(\la)\.\Ps_s\.z_s^{[j]}\.z_k^{[j]*}\\
     &\hspace*{35mm}+\sum_{s\in\Iz}\sum_{\la_j}\sum_{\la_i} (\la-\la_j)^{-1}\.(\bla-\la_i)^{-1}\.
              z_k^{[j]}\.z_s^{[j]*}\.\Ps_s\.z_s^{[i]}\.z_k^{[i]*}\\
    &\hspace*{5mm}
     =\sum_{s\in\Iz} G_{k,s}(\la)\.\Ps_s\.G_{k,s}^*(\la)
      -2\sum_{\la_j} (\la-\la_j)^{-1}\.(\bla-\la_j)^{-1}\.z_k^{[j]}\.z_k^{[j]*}\\
      &\hspace*{20mm}+\sum_{\la_j}\sum_{\la_i} \de_{ji}\.(\la-\la_j)^{-1}\.(\bla-\la_i)^{-1}\.z_k^{[j]}\.z_k^{[i]*}\\
    &\hspace*{5mm}
     =\sum_{s\in\Iz} G_{k,s}(\la)\.\Ps_s\.G_{k,s}^*(\la)-\sum_{\la_j} \abs{\la-\la_j}^{-2}\.z_k^{[j]}\.z_k^{[j]*},
  \end{align*}
 which implies~\eqref{E:EE14b}. The second inequality is a simple consequence of the latter inequality and the fact 
 that $\tr uv^*=\tr v^*u=v^*u$ for any vectors $u,v\in\Cbb^{2n}$.
\end{proof}


\section{Main results}\label{S:main}

If we take the finite orthonormal system of eigenfunctions $z^{[1]},z^{[2]},\dots,z^{[r]}$ for~\eqref{E:Sla} 
constructed in the preceding section and arbitrary $c_1^{[1]},\dots,c_1^{[r_1]},c_2^{[1]},\dots\in\Cbb$, then
 \begin{equation*}
  z\coloneq \sum_{j=1}^{r}\sum_{\ell=1}^{r_j} c_j^{[\ell]}\.z^{[\ell]}(\la_j)
 \end{equation*}
is a well-defined sequence from $\Cbb(\Izp)^{2n}$. The connection between $z$ and the coefficients $c_j$ can be easily 
derived as in the case of the Fourier series expansion, i.e., 
 \begin{equation}\label{E:coeff}
  \innerP{z^{[\ell]}(\la_j)}{z}
   =\sum_{k\in\Iz} z_k^{[\ell]*}(\la_j)\.\Ps_k\.\big(c_1^{[1]}\.z^{[1]}_k(\la_1)
                                                   +\dots+c_j^{[\ell]}\.z^{[\ell]}_k(\la_j)+\dots+\big)=c_j^{[\ell]},
 \end{equation}
which leads to the expression
 \begin{equation}\label{E:EE1}
  z\sim\sum_{j=1}^{r}\sum_{\ell=1}^{r_j} z^{[\ell]}(\la_j)\.\sum_{k\in\Iz} z_k^{[\ell]*}(\la_j)\.\Ps_k\.z_k.
 \end{equation}
However, can we replace the symbol $\sim$ by $=$? And can we get an arbitrary element of $\Cbb(\Izp)^{2n}$ on the 
right-hand side? The answer is, of course, negative because of several reasons. The first one is related to the 
definition of eigenfunctions of~\eqref{E:Sla}, i.e., the sequence obtained through this expression always satisfies the 
same boundary values $\al\.z_0=0=\be\.z_{N+1}$. The second reason is related to the absence of the index $N+1$ in the 
third sum on the right-hand side of~\eqref{E:EE1}, i.e., any change of the value of $z_{N+1}$ does not affect the value 
of $c_j^{[\ell]}$ and, consequently, the right-hand side of~\eqref{E:EE1}. In addition, all coefficients $c_j^{[\ell]}$
can be zero, while $z\not\equiv0$. This happens especially when $z\in[0]$. Actually, various representatives of $[z]$ 
yield the same right-hand side of~\eqref{E:EE1}. Therefore, the expression in~\eqref{E:EE1} should be rather written as
 \begin{equation}\label{E:EE2}
  [z]\sim\bigg[\sum_{j=1}^{r}\sum_{\ell=1}^{r_j} z^{[\ell]}(\la_j)\.
                                                    \sum_{k\in\Iz} z_k^{[\ell]*}(\la_j)\.\Ps_k\.z_k\bigg]
 \end{equation}
for the corresponding class $[z]\ni z$. Before we provide a complete description of classes $[z]\in\tltp$ being 
expressible as in~\eqref{E:EE2}, we need to prove the following lemmas.

\begin{lemma}\label{L:L8}
 Let $\al,\be\in\Ga$ be given and Hypothesis~\ref{H:WAC} hold. Furthermore, let $a>0$ and $f\in\Cbb(\Izp)^{2n}$ be 
 such that
  \begin{equation}\label{E:EE15}
   \sum_{k\in\Iz} z_k^{[\ell]*}(\la_j)\.\Ps_k\.f_k=0
  \end{equation}
 for all eigenfunctions corresponding to eigenvalues obeying $\abs{\la_j}\leq a$. If $\la_j\neq0$ for all 
 $j\in\{1,\dots,r\}$ and $\hz$ solves the nonhomogeneous problem~\eqref{E:EE4} with $\la=0$, i.e., it satisfies
  \begin{equation}\label{E:EE15a}
   \Slaf{0}{f},\quad k\in\Iz,\quad \al\.z_0=0=\be\.z_{N+1},
  \end{equation}
 then
  \begin{equation}\label{E:EE16}
   \normP{\hz}^2\leq a^{-2} \normP{f}^2.
  \end{equation}
 If $\la_j=0$ for some $j\in\{1,\dots,r\}$ and, in addition, the sequence $\hz$ satisfies
  \begin{equation}\label{E:EE16a}
   \sum_{k\in\Iz} z_k^{[\ell]*}(0)\.\Ps_k\.\hz_k=0
  \end{equation} 
 for all eigenfunctions $z_k^{[\ell]}(0)$ corresponding to the zero eigenvalue, then the latter estimate 
 in~\eqref{E:EE16} remains valid.
\end{lemma}
\begin{proof}
 We start with the first case when all eigenvalues are nonzero. Let $\la\in\Cbb$ not be an eigenvalue and 
 $y\in\Cbb(\Izp)^{2n}$ denote the unique solution of problem~\eqref{E:EE4}, i.e.,
  \begin{equation*}
   y_k(\la)=\sum_{j\in\Iz} G_{k,j}(\la)\.\Ps_j\.f_j,\quad k\in\Iz,
  \end{equation*}
 by Theorem~\ref{T:T5+C5A} and equality~\eqref{E:EE5a}. As the main tool of this proof we utilize the Taylor series 
 expansion of the function 
  \begin{equation*}
   \om(\la)\coloneq \sum_{k\in\Iz} \hz_k^*\.\Ps_k\.y_k(\la)
   =\sum_{k\in\Iz}\sum_{j\in\Iz} \hz_k^*\.\Ps_k\.G_{k,j}(\la)\.\Ps_j\.f_j
  \end{equation*}
 with $\hz=\sum_{j\in\Iz} G_{\cdot,j}(0)\.\Ps_j\.f_j$. We know from the previous section that the function $\om(\la)$ 
 is holomorphic in $\Cbb$ except for all eigenvalues, where it has a simple pole and the residue at the eigenvalue 
 $\la_s$ is, by Theorem~\ref{T:T7}, equal to
  \begin{align}
   \lim_{\la\to\la_s} (\la-\la_s)\.\om(\la)
    &=\sum_{k\in\Iz}\sum_{j\in\Iz} \hz_k^*\.\Ps_k\big[\lim_{\la\to\la_s}(\la-\la_s)\.G_{k,j}(\la)\big]\.\Ps_j\.f_j 
                                                                                                               \notag\\
    &=-\sum_{k\in\Iz}\sum_{j\in\Iz}
           \sum_{\ell=1}^{r_s} \hz_k^*\.\Ps_k\.z_k^{[\ell]}(\la_s)\.z_j^{[\ell]*}(\la_s)\.\Ps_j\.f_j,   \label{E:EE17}
  \end{align}
 where the third sum in the last expression includes all eigenfunctions from the orthonormal set corresponding to 
 $\la_s$. Since the right-hand side of~\eqref{E:EE17} is zero for any eigenvalue $\abs{\la_s}\leq a$ 
 by~\eqref{E:EE15}, it means that these eigenvalues are removable singularities of $\om(\la)$, i.e.,
 this function possesses a power series expansion in $\abs{\la}\leq a$ and it can be made holomorphic in this disk by 
 the Riemann theorem. We aim to derive this expansion in the next part of the proof. But before that we need to focus 
 on $y_k(\la)$ and its possible expression as
  \begin{equation}\label{E:EE18}
   y_k(\la)=\sum_{m=0}^\infty \la^m\.y_k^{[m]}
  \end{equation}
 for all $\abs{\la}\leq a$ small enough and all $k\in\Izp$. Upon inserting~\eqref{E:EE18} back into~\eqref{E:EE4} we 
 get the boundary value problem
  \begin{gather*}
   \sum_{m=0}^\infty \la^m\.y_k^{[m]}=\sum_{m=0}^\infty \la^m\.\Sc_k\.y_{k+1}^{[m]}
    -\Jc\.\Ps_k\Big(f_k+\sum_{m=0}^\infty \la^{m+1}\.\Sc_k\.y_{k+1}^{[m]}\Big),\\
   \al\.\sum_{m=0}^\infty \la^m\.y_0^{[m]}=0=\be\.\sum_{m=0}^\infty \la^m\.y_{N+1}^{[m]},
  \end{gather*}
 which means that the coefficients $y_k^{[0]},y_k^{[1]},\dots$ in~\eqref{E:EE18} can be obtained as the unique 
 solutions of the infinite systems of boundary value problems
  \begin{equation}\label{E:EE18a}
   \left. 
    \begin{aligned}
     &y_k^{[0]}=\Sc_k\.y_{k+1}^{[0]}-\Jc\.\Ps_k\.f_k,\quad k\in\Iz,\quad \al\.y_0^{[0]}=0=\be\.y_{N+1}^{[0]},\\
     &y_k^{[1]}=\Sc_k\.y_{k+1}^{[1]}-\Jc\.\Ps_k\.y_k^{[0]},\quad k\in\Iz,\quad \al\.y_0^{[1]}=0=\be\.y_{N+1}^{[1]},\\
     &y_k^{[2]}=\Sc_k\.y_{k+1}^{[2]}-\Jc\.\Ps_k\.y_k^{[1]},\quad k\in\Iz,\quad \al\.y_0^{[2]}=0=\be\.y_{N+1}^{[2]},\\
     &\hspace*{50mm}\vdots
    \end{aligned}
   \,\,\,\right\}
  \end{equation}
 with the nonhomogeneous parts given recursively, i.e., $y^{[0]}$ is the unique solution of system~\Slaf{0}{f}, while 
 $y^{[m]}$ solves~\Slaf{0}{y^{[m-1]}} for $m=1,2,\dots$ More precisely, 
  \begin{equation}\label{E:EE19}
   \left. 
    \begin{aligned}
     &y^{[0]}_k=\hz_k=\sum_{j\in\Iz} G_{k,j}(0)\.\Ps_j\.f_j,\\
     &y^{[1]}_k=\sum_{j\in\Iz} G_{k,j}(0)\.\Ps_j\.y^{[0]}
            =\sum_{j\in\Iz}\sum_{s\in\Iz} G_{k,j}(0)\.\Ps_j\.G_{j,s}(0)\.\Ps_s\.f_s,\\
     &y^{[2]}_k=\sum_{j\in\Iz} G_{k,j}(0)\.\Ps_j\.y^{[1]}
               =\sum_{j\in\Iz}\sum_{s\in\Iz}\sum_{p\in\Iz} 
                             G_{k,j}(0)\.\Ps_j\.G_{j,s}(0)\.\Ps_s\.G_{s,p}(0)\.\Ps_p\.f_p,\\       
     &\hspace*{50mm}\vdots                           
    \end{aligned}
   \,\,\,\right\}
  \end{equation}
 by Theorem~\ref{T:T5+C5A} and equality~\eqref{E:EE5a}. Since $\la=0$ is not an eigenvalue, the function $y_k(\cdot)$ 
 is holomorphic in a neighborhood of the origin, where it admits a Taylor series expansion. It justifies the expression 
 of $y_k(\la)$ as in~\eqref{E:EE18} for all $\la$ small enough with the coefficients specified in~\eqref{E:EE19}. 
 Consequently, we have
  \begin{equation*}
   \om(\la)=\sum_{m=0}^\infty \la^m\.\sum_{k\in\Iz} y_k^{[0]*}\.\Ps_k\.y_k^{[m]}
  \end{equation*}
 for all $\la$ small enough. Since, without loss of generality, $\om(\la)$ is holomorphic for $\abs{\la}\leq a$, this 
 expression is valid for all $\la$ in this disk. Therefore, by the Cauchy inequality for the Taylor series 
 coefficients, it holds
  \begin{equation}\label{E:EE20}
   \abs[\Big]{\sum_{k\in\Iz} y_k^{[0]*}\.\Ps_k\.y_k^{[m]}}\leq b\.a^{-m}, 
   \qtext{where $b\coloneq \max_{\abs{\la}=a}\abs{\om(\la)}$.}
  \end{equation}
 
 Now, if we put $y^{[-1]}\coloneq f$, then the extended Lagrange identity~\eqref{E:lagrange} yields
  \begin{equation*}
   y_k^{[m]}\.\Jc\.y_k^{[s]}\bigg|_0^{N+1}
    =\sum_{k\in\Iz} \big(y_k^{[m-1]*}\.\Ps_k\.y_k^{[s]}-y_k^{[m]*}\.\Ps_k\.y_k^{[s-1]}\big)
  \end{equation*}
 for any pair $m,s\in\{0,1,\dots\}$. But the left-hand side is zero due to the boundary conditions in~\eqref{E:EE18a}, 
 so it follows
  \begin{equation*}
   \sum_{k\in\Iz} y_k^{[m-1]*}\.\Ps_k\.y_k^{[s]}=\sum_{k\in\Iz} y_k^{[m]*}\.\Ps_k\.y_k^{[s-1]}
  \end{equation*}
 for all $m,s\in\{0,1,\dots\}$. In particular, for the choice $m=s-1$ we achieve
  \begin{equation}\label{E:EE22}
   \normP{y^{[m]}}^2
    =\sum_{k\in\Iz} y_k^{[m-1]*}\.\Ps_k\.y_k^{[m+1]}
    =\dots=\sum_{k\in\Iz} y_k^{[0]*}\.\Ps_k\.y_k^{[2m]},
  \end{equation}
 which together with~\eqref{E:EE20} implies 
  \begin{equation}\label{E:EE23}
   \normP{y^{[m]}}^2=\sum_{k\in\Iz} y_k^{[0]*}\.\Ps_k\.y^{[2m]}_k\leq b\.a^{-2m}.
  \end{equation}
 At the same time, the Cauchy--Schwarz inequality applied to the second term in~\eqref{E:EE22} yields
  \begin{equation}\label{E:EE24}
   \normP{y^{[m]}}^2
    =\innerP{y^{[m-1]}}{y^{[m+1]}}
    \leq \normP{y^{[m-1]}}\.\normP{y^{[m+1]}}
  \end{equation}
 for all $m\in\{0,1,\dots\}$. If $\normP{y^{[0]}}=\normP{\hz}=0$, then inequality~\eqref{E:EE16} stated in the lemma is 
 obviously satisfied. In the opposite case, all norms $\normP{y^{[-1]}},\normP{y^{[1]}},\dots$ are also positive by 
 induction applied to~\eqref{E:EE24}, so we may calculate the ratio
  \begin{equation*}
   \si^{[m]}\coloneq \frac{\normP{y^{[m]}}^2}{\normP{y^{[m-1]}}^2},\quad m\in\{0,1,\dots\}.
  \end{equation*}
 The value $\si^{[m]}$ is nondecreasing by~\eqref{E:EE24}, i.e., 
 $\si^{[0]}\leq \si^{[1]}\leq \dots\leq \si^{[m]}\leq\cdots$, from which we get
  \begin{equation*}
   \normP{y^{[1]}}^2\geq \si^{[0]}\.\normP{y^{[0]}}^2,\quad
   \normP{y^{[2]}}^2\geq \si^{[0]}\.\normP{y^{[1]}}^2\geq \big(\si^{[0]}\big)^2\.\normP{y^{[0]}}^2,\quad\dots
  \end{equation*}
 Thus, by using~\eqref{E:EE23}, we obtain
  \begin{equation*}
   (\si^{[0]})^m\.\normP{y^{[0]}}^2\leq \normP{y^{[m]}}^2\leq b\.a^{-2m}
   \qtextq{or equivalently} \si^{[0]}\.\normP{y^{[0]}}^{2/m}\leq b^{1/m}\.a^{-2}
  \end{equation*}
 for all $m\in\{1,2,\dots\}$, which for $m\to\infty$ leads to $\normP{\hz}^2/\normP{f}^2=\si^{[0]}\leq a^{-2}$ as 
 stated in~\eqref{E:EE16}. 

 If $\la=0$ is an eigenvalue of~\eqref{E:Sla}, we modify the eigenvalue problem so that the new problem possesses only 
 nonzero eigenvalues and, consequently, the first part of the statement can be utilized. More precisely, we increase 
 all eigenvalues by a suitable small $\eps>0$, i.e., we replace system~\eqref{E:Sla} by
  \begin{equation}\label{E:hSla}\tag{$\hat{\text{S}}_\la$}
   z_k(\la)=(\hSc_k+\hla\.\hVc_k)\.z_{k+1}(\la)
  \end{equation}
 with $\hSc_k\coloneq \Sc_k-\eps\.\Vc_k$, $\hVc_k\coloneq \Vc_k$, and $\hla=\la+\eps$. The reader may easily verify 
 that this perturbed system meets all the basic assumptions as the original system. If $\eps\neq-\la_j$ for 
 all eigenvalues such that $\la_j<0$, then eigenvalues of system~\eqref{E:hSla} satisfy $\hla_j=\la_j+\eps\neq0$ and 
 their eigenspaces coincide with $\Eigen(\la_j)$, $j\in\{1,\dots,r\}$. The sequence $\hz$ now solves 
 system~\hSlaf{0}{\hf} with $\hf\coloneq f+\eps\.\hz$ and, in order to apply the first part of the statement, we need 
 to prove that condition~\eqref{E:EE15} with $f$ replaced by $\hf$ remains true, i.e., we need to show that 
  \begin{equation*}
   \sum_{k\in\Iz} z_k^{[\ell]*}(\hla_j)\.\Ps_k\.(f_k+\eps\.\hz_k)=0
  \end{equation*}
 for all eigenfunctions corresponding to eigenvalues $\abs{\hla_j}\leq \hat{a}$ with a suitable $\hat{a}>0$. 
 Actually, if we take $\hat{a}\coloneq a-\eps>0$, then according to the original assumption in~\eqref{E:EE15} and the 
 additional condition in~\eqref{E:EE16a} it suffices to verify only that
  \begin{equation*}
   \sum_{k\in\Iz} z_k^{[\ell]*}(\hla_j)\.\Ps_k\.\hz_k=0 
   \qtext{for all $\abs{\hla_j}\leq \hat{a}$ and $\hla_j\neq\eps$.}
  \end{equation*}
 However, this follows from~\eqref{E:EE15} and the extended Lagrange identity~\eqref{E:lagrange}, because
  \begin{align*}
   0=z_k^{[\ell]*}(\hla_j)\.\Jc\.\hz_k\Big|_{0}^{N+1}
    &=\sum_{k\in\Iz}\Big\{\hla_j\.z_k^{[\ell]*}(\hla_j)\.\Ps_k\.\hz_k-z_k^{[\ell]*}(\hla_j)\.\Ps_k\.\hf_k\Big\}\\
    &=\sum_{k\in\Iz}\Big\{(\la_j+\eps)\.z_k^{[\ell]*}(\hla_j)\.\Ps_k\.\hz_k-z_k^{[\ell]*}(\hla_j)\.\Ps_k\.f_k
      -\eps\.z_k^{[\ell]*}(\hla_j)\.\Ps_k\.\hz_k\Big\}
  \end{align*}
 or equivalently
  \begin{equation*}
   \la_j\.\sum_{k\in\Iz} z_k^{[\ell]*}(\hla_j)\.\Ps_k\.\hz_k=\sum_{k\in\Iz} z_k^{[\ell]*}(\hla_j)\.\Ps_k\.f_k
  \end{equation*}
 with $\la_j\neq0$ and the right-hand side being equal to $\sum_{k\in\Iz} z_k^{[\ell]*}(\la_j)\.\Ps_k\.f_k$ and, 
 consequently, zero by~\eqref{E:EE15}. Therefore, all the assumptions of the first part of the lemma applied to 
 system~\eqref{E:hSla} are satisfied, and so we have the estimate
  \begin{equation*}
   \normP{\hz}^2\leq \hat{a}^{-2}\.\normP{\hf}^2.
  \end{equation*}
 Since this inequality is true for any $\eps>0$ small enough, by letting $\eps\to0^+$ we obtain (again) 
 inequality~\eqref{E:EE16}, which completes the proof.
\end{proof}

As a direct consequence of Lemma~\ref{L:L8} we obtain the following statement.

\begin{lemma}\label{T:T9}
 Let $\al,\be\in\Ga$ be given, Hypothesis~\ref{H:WAC} hold, and $f\in\Cbb(\Izp)^{2n}$ be arbitrary. 
 If $\hz\in\Cbb(\Izp)^{2n}$ is a solution of~\eqref{E:EE15a}, then for any $a>0$ the sequence $\hz^a\in\Cbb(\Izp)^{2n}$ 
 given by
  \begin{equation*}
   \hz_k^a\coloneq \hz_k-\sum_{\abs{\la_j}\leq a}\sum_{\ell=1}^{r_j} c_j^{[\ell]}\.z_k^{[\ell]}(\la_j),\quad k\in\Izp,
  \end{equation*}
 with the Fourier coefficients as in~\eqref{E:coeff}, i.e.,
  \begin{equation*}
   c_j^{[\ell]}\coloneq \sum_{k\in\Iz} z_k^{[\ell]*}(\la_j)\.\Ps_k\.\hz_k,
  \end{equation*}
 satisfies the inequality
  \begin{equation}\label{E:EE25}
   \normP{\hz^{a}}^2\leq a^{-2} \normP{f}^2.
  \end{equation}
\end{lemma}
\begin{proof}
 For all eigenvalues satisfying $\abs{\la_j}\leq a$ and $\ell=1,\dots,r_j$ we define the numbers
  \begin{equation*}
   d_j^{[\ell]}\coloneq \sum_{k\in\Iz} z_k^{[\ell]*}(\la_j)\.\Ps_k\.f_k,
  \end{equation*}
 which determine the sequence
  \begin{equation*}
   f_k^a\coloneq f_k-\sum_{\abs{\la_j}\leq a}\sum_{\ell=1}^{r_j} d_j^{[\ell]}\.z_k^{[\ell]}(\la_j), \quad k\in\Izp.
  \end{equation*}
 Since $z^{[\ell]}(\la_j)$ and $\hz$ satisfy the same boundary conditions, the extended Lagrange 
 identity~\eqref{E:lagrange} yields
  \begin{equation*}
   0=z_k^{[\ell]*}(\la_j)\.\Jc\.\hz_k\Big|_0^ { N+1 } 
    =\sum_{k\in\Iz} \big\{\la_j\.z_k^{[\ell]*}(\la_j)\.\Ps_j\.\hz_k-z_k^{[\ell]*}(\la_j)\.\Ps_k\.f_k\big\}
  \end{equation*}
 or equivalently $\la_j\.c_j^{[\ell]}=d_j^{[\ell]}$. Consequently, it follows 
  \begin{equation*}
   \hz_k^a=\Sc_k\Big(\hz_{k+1}-\sum_{\abs{\la_j}\leq a}\sum_{\ell=1}^{r_j} c_j^{[\ell]}\.z_{k+1}^{[\ell]}(\la_j)\Big)
            -\Jc\.\Ps_k\Big(f_k-\sum_{\abs{\la_j}\leq a}\sum_{\ell=1}^{r_j} 
                                                \la_j\.c_j^{[\ell]}\.\Sc_k\.z_{k+1}^{[\ell]}(\la_j)\Big)
  \end{equation*}
 i.e., the sequence $\hz^a$ solves system~\Slaf{0}{f^a}, and simultaneously
  \begin{gather*}
   \al\.\hz_0^a=\al\.\hz_0-\sum_{\abs{\la_j}\leq a}\sum_{\ell=1}^{r_j} c_j^{[\ell]}\.\al\.z_0^{[\ell]}(\la_j)=0 
   \intertext{and}
   \be\.\hz_{N+1}^a=\be\.\hz_{N+1}-
                         \sum_{\abs{\la_j}\leq a}\sum_{\ell=1}^{r_j} c_j^{[\ell]}\.\be\.z_{N+1}^{[\ell]}(\la_j)=0.
  \end{gather*}
 Furthermore, the orthonormality of the set of eigenfunctions implies
  \begin{equation*}
   \sum_{k\in\Iz} z_k^{[\ell]*}(\la_j)\.\Ps_k\.\hz_k^a=c_j^{[\ell]}-c_j^{[\ell]}=0 \qtextq{and}
   \sum_{k\in\Iz} z_k^{[\ell]*}(\la_j)\.\Ps_k\.f_k^a=d_j^{[\ell]}-d_j^{[\ell]}=0.
  \end{equation*}
 Thus, the sequence $\hz^a$ meets all the assumptions of Lemma~\ref{L:L8} (whether $\la=0$ is an eigenvalue or 
 not) with $f$ replaced by $f^a$, which yields
  \begin{align*}
   \normP{\hz^{a}}^2
    &\leq a^{-2} \normP{f^{a}}^2\\
    &=a^{-2} \Big\{\normP{f}^2
      -\sum_{\abs{\la_j}\leq a}\sum_{\ell=1}^{r_j}\sum_{k\in\Iz} 
                                      d_j^{[\ell]*}z_k^{[\ell]*}(\la_j)\.\Ps_k\.f_k\\
     &\hspace*{35mm}-\sum_{\abs{\la_j}\leq a}\sum_{\ell=1}^{r_j}\sum_{k\in\Iz}
                                     d_j^{[\ell]}\.\Big(z_k^{[\ell]*}(\la_j)\.\Ps_k\.f_k\Big)^{\!\!*}\\
     &\hspace*{45mm}+\sum_{\abs{\la_j}\leq a}\sum_{\ell=1}^{r_j}\sum_{\abs{\la_m}\leq a}\sum_{s=1}^{r_m}
                     \sum_{k\in\Iz}\bar{d_j^{[\ell]}}\.d_r^{[s]}\.z_k^{[\ell]*}(\la_j)\.\Ps_k\.z_k^{[s]}(\la_r)\Big\}\\
   &\hspace*{55mm}=a^{-2}\.\Big(\normP{f}^2
   -\sum_{\abs{\la_j}\leq a}\sum_{\ell=1}^{r_j}\abs[\big]{d_j^{[\ell]}}^2\Big) \leq a^{-2}\normP{f}^2
  \end{align*}
 as stated in~\eqref{E:EE25}.
\end{proof}

As we let $a\to\infty$ in Lemma~\ref{T:T9}, the right-hand side of~\eqref{E:EE25} tends to zero for any 
$f\in\Cbb(\Izp)^{2n}$, which means that $\hz$ and $\sum_{j=1}^r\sum_{\ell=1}^{r_j} c_j^{[\ell]}\.z_k^{[\ell]}(\la_j)$
belong to the same equivalence class. This conclusion is our first main result, which is formulated below together with 
Parseval's identity as a simple consequence of equality~\eqref{E:expansion.norm}. Of course, this 
``mean square convergence'' does not imply that these two sequences are the same as we illustrate in 
Example~\ref{Ex:expansion} below, see also Corollary~\ref{C:expansion}. This equality cannot be forced by 
the choice $\Ps_k\equiv I$ as in the last part of Theorem~\ref{T:EE.admissible}, because it contradicts 
Hypothesis~\ref{H:basic}. Similarly, the conclusion does not change if we replace the Weak Atkinson condition by the 
Strong Atkinson condition, because $\hz$ solves the nonhomogeneous problem. Instead of that we should employ much more 
restrictive alternative of the Atkinson-type condition dealing with system~\Slaf{0}{f}, which would be analogous to the 
condition used in~\cite[Theorem~7.6]{amK89:SJMA.a} for system~\eqref{E:cLHS}. Consequently, it would imply the 
existence 
of a~densely defined operator associated with the mapping $\mL(z)$ introduced in~\eqref{E:equiv.Sla}. However, this 
condition is never satisfied in the setting of discrete symplectic systems except for the most trivial case 
$\Ps_k\equiv0$ as it was discussed in~\cite{pZ20}.


\begin{theorem}\label{T:main.expansion}
 Let $\al,\be\in\Ga$ be given, Hypothesis~\ref{H:WAC} hold, and $f\in\Cbb(\Izp)^{2n}$ be arbitrary. Then the
 equivalence class $[\hz]\in\tltp$ determined by a solution $\hz\in\Cbb(\Izp)^{2n}$ of problem~\eqref{E:EE15a} 
 possesses the representative given as the sum
  \begin{equation}\label{E:expansion.coeff}
   \sum_{j=1}^r\sum_{\ell=1}^{r_j} c_j^{[\ell]}\.z^{[\ell]}(\la_j)
   \qtextq{with the coefficients} c_j^{[\ell]}\coloneq \sum_{k\in\Iz} z_k^{[\ell]*}(\la_j)\.\Ps_k\.\hz_k,
  \end{equation}
 where $z^{[1]}(\la_1),\dots,z^{[r_1]}(\la_1),z^{[1]}(\la_2),\dots\in\Cbb(\Izp)^{2n}$ form the finite orthonormal set 
 of eigenfunctions of system~\eqref{E:Sla} established in Theorem~\ref{T:orthonormal.set}, i.e., the solution $\hz$ 
satisfies 
  \begin{equation}\label{E:expansion.norm}
   \normP[\bigg]{\hz- \sum_{j=1}^r \sum_{\ell=1}^{r_j} c_j^{[\ell]}\.z^{[\ell]}(\la_j)}=0
  \end{equation}
 and, consequently, its norm and the coefficients are related through Parserval's identity 
  \begin{equation}\label{E:parseval}
   \normP{\hz}^2=\sum_{j=1}^r\sum_{\ell=1}^{r_j} \abs[\big]{c_j^{[\ell]}}^2.
  \end{equation} 
\end{theorem}

\begin{remark}\label{R:EE.theorems.compare}\leavevmode
 Can we get Theorem~\ref{T:EE.admissible} as a special case of Theorem~\ref{T:main.expansion}? 
 System~\eqref{E:dss.blocks} is equivalent to~\eqref{E:Sla} with $\Sc_k\coloneq-\Jc\.\tSc_k^*\.\Jc$ and 
 $\Vc_k\coloneq-\Jc\.\tSc_k^*\.\diag\{\tWc_k,0\}$ or to its alternative form displayed in~\eqref{E:equiv.Sla} with the 
 weight matrix
  \begin{equation*}
   \Ps_k\coloneq\mmatrix{\tAc_k^*\\[1mm] \tBc_k^*}\.\tWc_k\.\mmatrix{\tAc_k & \tBc_k},
  \end{equation*}
 which implies the equality $z_{k+1}^*(\la)\.\diag\{\tWc_k,0\}\.z_{k+1}(\la)=z_k^*(\la)\.\Ps_k\.z_k(\la)$ for every 
 admissible sequence for~\eqref{E:dss.blocks}. Furthermore, an admissible sequence 
 $\hz=(\hx^*,\hu^*)^*$ for~\eqref{E:dss.blocks} satisfying $\hx_0=0=\hx_{N+1}$ can be seen as a 
 solution of problem~\eqref{E:EE15a} with the coefficient matrices specified above and $\al=(I,0)=\be$ if and only if 
 there exists a sequence $f\in\Cbb(\Izp)^{2n}$ such that 
  \begin{equation*}
   u_{k+1}-\tCc_k\.x_k-\tDc_k\.u_k=\mmatrix{\tWc_k\.\tAc_k & \tWc_k\.\tBc_k}\!f \qtext{for all $k\in\Iz$,}
  \end{equation*}
 i.e., if and only if $u_{k+1}-\tCc_k\.x_k-\tDc_k\.u_k$ belongs to the range of the matrix 
 $\tWc_k\.(\tAc_k\,\,\,\,\tBc_k)$ for all $k\in\Iz$. It shows that Theorems~\ref{T:EE.admissible} 
 and~\ref{T:main.expansion} coincide in the case $\tWc_k>0$ on $\Iz$, because $\rank(\tAc_k\,\,\,\,\tBc_k)\equiv n$ on 
 $\Iz$. Otherwise, Theorem~\ref{T:EE.admissible} applies to a larger class of sequences than we can get through this 
 special form of problem~\eqref{E:EE15a}.

\end{remark}

We illustrate the conclusion of Theorem~\ref{T:main.expansion} in the following example dealing with the same systems 
as in Example~\ref{Ex:eigenvalue}.

\begin{example}\label{Ex:expansion}\leavevmode
 \begin{enumerate}[leftmargin=10mm,topsep=0mm,label={{\normalfont{(\roman*)}}}]
  \item Let us start with system~\eqref{E:Sla} on an arbitrary finite discrete interval $\Iz=\onZ$ with the coefficient 
        matrices
         \begin{equation*}
          \Sc_k\equiv I_2\in\Cbb^{2\times2} \qtextq{and} \Ps_k=\diag\{0,\De v_k\},
         \end{equation*}
        where $v\in\Cbb(\Izp)^1$ is a given non-decreasing real-valued scalar sequence with $v_0\coloneq0$ and 
        $v_k>0$ for some $k\in\Izp\stm\{0\}$, which guarantees that Hypothesis~\ref{H:WAC} holds. If we take 
        $\al=(1\,\,\,\,0)=\be$ then $\la_1=0$ is the only eigenvalue and the corresponding normalized eigenfunction is
         \begin{equation*}
          z_k^{[1]}(\la_1)\equiv\mmatrix{0\\ 1/\sqrt{v_{N+1}}},
         \end{equation*}
        i.e., $\eta_1^{[1]}=1/\sqrt{v_{N+1}}$. For an arbitrary sequence 
        $f=\msmatrix{f^{[1]}\\ f^{[2]}}\in\Cbb(\Izp)^2$ problem~\eqref{E:EE15a} leads to
         \begin{equation}\label{E:example.expansion.i}
          \De x_k=f_k^{[2]}\.\De v_k \qtextq{and}
          \De u_k=0,\quad k\in\Iz,\quad x_0=0=x_{N+1}.
         \end{equation}
        Since $\la=0$ is an eigenvalue, the assumptions of Theorem~\ref{T:T5+C5A} are not satisfied, and so the 
        existence of a solution of~\eqref{E:example.expansion.i} may be violated. Actually, a solution exists if and 
        only if $\sum_{k\in\Iz} f_k^{[2]}\.\De v_k=0$, in which case we have
         \begin{equation*}
          \hz_k=\mmatrix{\sum_{j=0}^{k-1} f_j^{[2]}\.\De v_j \\ u},\quad k\in\Izp,
         \end{equation*}
        for an arbitrary constant $u\in\Cbb$. Then the unique ``Fourier'' coefficient from~\eqref{E:expansion.coeff} 
        is equal to $c_1^{[1]}=\sum_{k\in\Iz} z_k^{[1]}(0)\.\Ps_k\.\hz_k=u\.\sqrt{v_{N+1}}$, which yields the 
        eigenfunctions expansion as the constant sequence 
         \begin{equation*}
          c_1^{[1]}\.z_k^{[1]}(0)\equiv\mmatrix{0\\ u}, \quad k\in\Izp.
         \end{equation*}
        Consequently, $\hz-c_1^{[1]}\.z_k^{[1]}(0)=\msmatrix{\sum_{j=0}^{k-1} f_j^{[2]}\.\De v_j\\ 0}\in[0]$ and
        $\normP{\hz}^2=\abs[\big]{c_1^{[1]}}^2=\abs{u}^2\.v_{N+1}$ as stated in Theorem~\ref{T:main.expansion}.
        
        On the other hand, for $\al=(1\,\,\,\,0)$ and $\be=(0\,\,\,\,1)$ there is no eigenvalue, while a solution of 
        problem~\eqref{E:EE15a} reducing to
         \begin{equation*}
          \De x_k=f_k^{[2]}\.\De v_k \qtextq{and}
          \De u_k=0,\quad k\in\Iz,\quad x_0=0=u_{N+1},
         \end{equation*}
        exists for any $f\in\Cbb(\Izp)^2$ with $\hz_k=\msmatrix{\sum_{j=0}^{k-1} f_j^{[2]}\.\De v_j \\ 0}$, $k\in\Izp$.
        Consequently, $\hz\in[0]$, which agrees with the conclusion of Theorem~\ref{T:main.expansion}.
        
  \item Now let us take an arbitrary finite discrete interval $\Iz=\onZ$ and system~\eqref{E:Sla} with
         \begin{equation*}
          \Sc_k\equiv I_4\in\Cbb^{4\times4} \qtextq{and} 
          \Ps_k\equiv\mmatrix{a\.I_2 & \sqrt{ab}\.I_2\\ \sqrt{ab}\.I_2 & b\.I_2},
         \end{equation*}
        where $a,b>0$ are given real numbers, which guarantees that Hypothesis~\ref{H:WAC} holds. If we take
        $\al=(I_2\,\,\,\,0)$ and $\be=\msmatrix{1 & 0 & 0 & 0\\ 0 & 0 & 0 & 1}$, then we have two simple eigenvalues
        $\la_1=0$ and $\la_2=1/[(N+1)\.\sqrt{ab}]$ with the corresponding orthonormal eigenfunctions
         \begin{equation*}
          z_k(\la_1)\equiv \mmatrix{0\\ 0\\ \frac{1}{\sqrt{b\.(N+1)}}\\ 0}\qtextq{and}
          z_k(\la_2)=\mmatrix{0\\ \frac{k}{\sqrt{a\.(N+1)^3}}\\ 0\\
                                                  \frac{1}{\sqrt{b\.(N+1)}}-\frac{k}{\sqrt{b\.(N+1)^3}}},
          \qtext{$k\in\Izp$.} 
         \end{equation*}
        Problem~\eqref{E:EE15a} reduces to
         \begin{equation*}
          \De z_k=\mmatrix{\sqrt{ab}\.I_2 & b\.I_2\\ -a\.I_2 & -\sqrt{ab}\.I_2}\.\mmatrix{f_k^{[1]}\\ f_k^{[2]}},\quad 
          k\in\Iz,\quad x_0=0=x^{[1]}_{N+1}=u^{[2]}_{N+1}
         \end{equation*}
        with the notation
        $f=\msmatrix{f^{[1]}\\ f^{[2]}}=\big(f^{[1,1]},f^{[1,2]},f^{[2,1]},f^{[2,2]}\big)^{\!\top}\in\Cbb(\Izp)^{4}$, 
        $z=\msmatrix{x\\ u}\in\Cbb(\Izp)^4$, 
        $x=\big(x^{[1]}, x^{[2]}\big)^{\!\top}\in\Cbb(\Izp)^2$, and 
        $u=\big(u^{[1]}, u^{[2]}\big)^{\!\top}\in\Cbb(\Izp)^2$. Its solution exists if and only if
        $\sqrt{ab}\.\sum_{k\in\Iz} f_k^{[1,1]}+b\.\sum_{k\in\Iz} f_k^{[2,1]}=0$, in which case it holds
         \begin{equation*}
          \hx_k=\sum_{j=0}^{k-1} \big(\sqrt{ab}\.f_j^{[1]}+b\.f_j^{[2]}\big) \qtextq{and}
          \hu_k=\hu_0-\sum_{j=0}^{k-1} \big(a\.f_j^{[1]}+\sqrt{ab}\.f_j^{[2]}\big),\quad k\in\Izp,
         \end{equation*}
        for an arbitrary $\hu_0^{[1]}\in\Cbb$ and 
        $\hu_0^{[2]}=\sum_{k\in\Iz} \big(a\.f_j^{[1,2]}+\sqrt{ab}\.f_j^{[2,2]}\big)$. In particular, if we choose
        $f_k^{[1]}=\msmatrix{k\\ k^2}$ and $f_k^{[1]}=\msmatrix{-\sqrt{a/b}\.k\\ 1}$ and $\hu_0^{[1]}=1$, then
        $\hu_0^{[2]}=(N+1)\.(2a\.N^2+a\.N+6\sqrt{ab})$, which yields
         \begin{equation*}
          \hz_k=\mmatrix{0\\ \sqrt{ab}\.(2k^3-3k^2+k)/6+kb\\ 1\\ 
                         (N-k+1)\.\big[3\.\sqrt{ab}+a\.N^2+a\.N\.(2k+1)/2+a\.k^2-a\.k/2\big]/3}, \quad k\in\Izp.
         \end{equation*}
        The corresponding ``Fourier'' coefficients are
         \begin{equation*}
          c_1^{[1]}=\sqrt{b\.(N+1)} \qtextq{and} 
          c_2^{[1]}=\sqrt{a\.b\.(N+1)^3}\.\big[6\.\sqrt{b}+\sqrt{a}\.N\.(2N+1)]/6,
         \end{equation*}
        and so we have the expansion
         \begin{equation*}
          c_1^{[1]}\.z_k^{[1]}(\la_1)+c_2^{[1]}\.z_k^{[1]}(\la_2)
          =\mmatrix{0\\ k\.\sqrt{ab}\.N\.(2N+1)/6+kb\\ 1\\ (N+1-k)\.\big[\sqrt{ab}+a\.N\.(2N+1)/6\big]}, \quad k\in\Izp.
         \end{equation*}
        Then, by a direct calculation, it can be verified that
         \begin{equation*}
          \hz_k-c_1^{[1]}\.z_k^{[1]}(\la_1)-c_2^{[1]}\.z_k^{[1]}(\la_2)
          =\mmatrix{0\\-\sqrt{ab}\.(N+1-k)\.(2N-1+2k)\.k/6 \\ 1\\ a\.k\.(2N^2-2k^2+N+3k-1)/6}\in[0]
         \end{equation*}
        and also
         \begin{align*}
          \normP{\hz}^2
           &=\abs[\big]{c_1^{[1]}}^2+\abs[\big]{c_2^{[1]}}^2\\
           &=\frac{(N+1)\.N\.ab}{36}\.\Big(4\.a\.N^5+24\.N^3\sqrt{a\.b}
               +12\.a\.{N}^4+60\.N^2\sqrt{a\.b}+13\.a\.N^3\\
           &\hspace{42mm}+48\.N\sqrt{a\.b}+6\.a\.N^2+12\.\sqrt{a\.b}+N\.a\\
           &\hspace{52mm}+36\.N\.b+72\.b\Big) +(ab+1)\.(N+1)\.b
         \end{align*}
        as stated in Theorem~\ref{T:main.expansion}.
 \end{enumerate}
\end{example}

Albeit Theorem~\ref{T:main.expansion} cannot provide a general ``pointwise'' eigenfunctions expansion due to the 
reasons mentioned above, the situation is not entirely hopeless. It is possible to achieve this type of expansion at 
least for some classes in $\tltp$ specified in the following corollary, compare with~\cite[Problem~9.16]{fvA64}. This 
result is a key ingredient in the proof of an integral representation of the $M(\la)$-function investigated in the 
subsequent part.

\begin{corollary}\label{C:expansion}
 Let $\al,\be\in\Ga$ be given, Hypothesis~\ref{H:WAC} hold, and $g\in\Cbb(\Izp)^{2n}$ be such that the equivalence 
 class $[g]$ possesses a representative being a solution of~\eqref{E:EE15a} for some $f\in\Cbb(\Izp)^{2n}$. Whenever 
 $\la\in\Cbb$ is not an eigenvalue of~\eqref{E:Sla}, then the unique solution $\tz\in\Cbb(\Izp)^{2n}$ of the boundary 
 value problem
  \begin{equation}\label{E:expansion.corollary.system}
   \Slaf{\la}{g},\quad k\in\Iz, \quad \al\.z_0=0=\be\.z_{N+1},
  \end{equation}
 admits the representation
  \begin{equation}\label{E:expansion.corollary}
   \tz_k=\sum_{j\in\Iz} G_{k,j}(\la)\.\Ps_j\.g_j
        =\sum_{j=1}^r \sum_{\ell=1}^{r_j} (\la_j-\la)^{-1}\. d_j^{[\ell]}\.z_k^{[\ell]}(\la_j) 
   \qtext{for all $k\in\Izp$}
  \end{equation}
 with the same eigenfunctions $z^{[1]}(\la_1),\dots,z^{[r_1]}(\la_1),z^{[1]}(\la_2),\dots\in\Cbb(\Izp)^{2n}$ as in 
 Theorem~\ref{T:main.expansion} and the coefficients given by 
  \begin{equation*}
   d_j^{[\ell]}\coloneq \sum_{k\in\Iz} z_k^{[\ell]*}(\la_j)\.\Ps_k\.g_k.
  \end{equation*}
\end{corollary}
\begin{proof}
 The existence and uniqueness of $\tz$ follow from Theorem~\ref{T:T5+C5A} as well as the first equality 
 in~\eqref{E:expansion.corollary}. To prove the second equality we need to show that the sequence on the right-hand 
 side solves problem~\eqref{E:expansion.corollary.system}. The validity of the boundary conditions is a 
 simple consequence of the definition of eigenfunctions of~\eqref{E:Sla}, while upon inserting this expression 
 into~\Slaf{\la}{g} we obtain
  \begin{equation}\label{E:expansion.corollary.proof}
   \begin{aligned}
    \tz_k-(\Sc_k+\la\.\Vc_k)\.\tz_{k+1}
     &=\sum_{j=1}^r \sum_{\ell=1}^{r_j} d_j^{[\ell]}\.\Vc_k\.z_{k+1}^{[\ell]}(\la_j)\\
     &=-\Jc\.\Ps_k\sum_{j=1}^r \sum_{\ell=1}^{r_j} d_j^{[\ell]}\.z_k^{[\ell]}(\la_j)
      =-\Jc\.\Ps_k\.g_k,
   \end{aligned}
  \end{equation}
 because the sequences $g$ and $\sum_{j=1}^r \sum_{\ell=1}^{r_j} d_j^{[\ell]}\.z_k^{[\ell]}(\la_j)$ belong to the same 
 equivalence class by Theorem~\ref{T:main.expansion}. More precisely, the last equality 
 in~\eqref{E:expansion.corollary.proof} is justified by the existence of $\hz\in[g]$ being a solution of~\eqref{E:EE15a}
 for some $f\in\Cbb(\Izp)^{2n}$, because in that case we have
 $\Ps_k\.g_k\equiv \Ps_k\.\hz_k\equiv \Ps_k\.\sum_{j=1}^r \sum_{\ell=1}^{r_j} c_j^{[\ell]}\.z_k^{[\ell]}(\la_j)$ on 
 $\Iz$ and $c_j^{[\ell]}=d_j^{[\ell]}$ by Theorem~\ref{T:main.expansion}.
\end{proof}

The right-hand side of the expansion given in~\eqref{E:EE1} can be alternatively written by using the {\it spectral} 
$n\times n$ matrix-valued {\it function}
 \begin{equation*}
  \tau_{\al,\be}(t)=\tau_{\al,\be,N}(t)\coloneq 
   \begin{cases}
    \sum_{0<\la_j\leq t}\sum_{\ell=1}^{r_j} \eta_j^{[\ell]}\.\eta_j^{[\ell]*}, & t>0,\\[1mm]
    -\sum_{t<\la_j\leq0}\sum_{\ell=1}^{r_j} \eta_j^{[\ell]}\.\eta_j^{[\ell]*}, & t<0,\\[1mm]
    0, & t=0,
   \end{cases}
 \end{equation*}
where $\eta_j^{[\ell]}\in\ker\be\.\tZ_{N+1}(\la_j)\in\Cbb^n$ are the vectors from the orthonormalization process 
discussed before Theorem~\ref{T:orthonormal.set}. The function $\tau_{\al,\be}(t)$ is a right-hand side continuous step 
function with jumps occurring only at the eigenvalues, it possesses the Hermitian property
 \begin{equation*}
  \tau_{\al,\be}^*(t)=\tau_{\al,\be}(t),
 \end{equation*}
and the nondecreasing property
 \begin{equation*}
  \tau_{\al,\be}(t)\leq \tau_{\al,\be}(s) \qtext{for any $t\leq s$.}
 \end{equation*}
Note that the value of the jump of $\tau_{\al,\be}(\cdot)$ at an eigenvalue $\la_j$ is equal to 
$\sum_{\ell=1}^{r_j} \eta_j^{[\ell]}\.\eta_j^{[\ell]*}$, which coincides (up to the sign) with the value of the residue 
$L_{-1}^{[j]}$ of $M_{N+1}(\la)$ at $\la_j$ calculated in Theorem~\ref{T:T7}. Since all eigenfunctions were constructed 
so that $z^{[\ell]}(\la_j)=\tZ(\la_j)\.\eta_j^{[\ell]}$, upon employing the spectral function and the 
Riemann--Stieltjes 
integral in the expansion given in~\eqref{E:EE1} we get
 \begin{align*}
  z
   &\sim\sum_{j=1}^r \sum_{\ell=1}^{r_j} 
              \tZ(\la_j)\.\eta_j^{[\ell]}\.\eta_j^{[\ell]*}\sum_{k\in\Iz}\tZ^*_k(\la_j)\.\Ps_k\.z_k\\
   &=\sum_{j=1}^r \sum_{\ell=1}^{r_j} 
              \tZ(\la_j)\.\big(\tau_{\al,\be}(\la_j)-\tau_{\al,\be}(\la_{j-1})\big)\om(\la_j)
   =\int_{-\infty}^\infty \tZ(\la)\. \d\tau_{\al,\be}(\la)\.\om(\la),
 \end{align*}
where we utilized Proposition~\ref{P:RS-integral} and the $n$-vector valued function
 \begin{equation*}\label{E:om.function.def}
  \om(\la)\coloneq \sum_{k\in\Iz}\tZ^*_k(\la)\.\Ps_k\.z_k.
 \end{equation*}
Furthermore, combining~\eqref{E:EE5b} and~\eqref{E:EE14b} with $k=0$ and \eqref{E:EE8} for $\im\la\neq0$ we obtain
 \begin{align}
  \sum_{j=1}^{r}\sum_{\ell=1}^{r_j} \abs{\la-\la_j}^{-2} z_0^{[\ell]}(\la_j)\.z_0^{[\ell]*}(\la_j)
   &\leq \sum_{s\in\Iz} G_{0,s}(\bla)\.\Ps_s\.G^*_{0,s}(\bla)\notag\\
   &=\sum_{s\in\Iz} \tZ_0(\bla)\.\Xc_s^*(\la)\.\Ps_s\.\Xc_s(\la)\.\tZ_0^*(\bla)\notag\\
   &=-\Jc\al^*\Big[\sum_{s\in\Iz} \Xc_s^*(\la)\.\Ps_s\.\Xc_s(\la)\Big]\.\al\.\Jc\notag\\
   &=-(\im\la)^{-1}\.\Jc\al^*\im M_{N+1}(\la)\.\al\.\Jc. \label{E:Im.Mla.estimate}
 \end{align}
Since $\abs{\la-\la_j}^{-2}\.\im\la=\im(\la_j-\la)^{-1}$, Proposition~\ref{P:RS-integral} implies that the left-hand 
side of~\eqref{E:Im.Mla.estimate} can be expressed as the Riemann--Stieltjes integral
 \begin{align}
  &\sum_{j=1}^{r}\sum_{\ell=1}^{r_j} \abs{\la-\la_j}^{-2} z_0^{[\ell]}(\la_j)\.z_0^{[\ell]*}(\la_j)\notag\\
   &\hspace*{20mm}=-(\im\la)^{-1}\.\Jc\al^*
     \bigg[\sum_{j=1}^{r}\sum_{\ell=1}^{r_j} \im(\la_j-\la)^{-1} \eta_j^{[\ell]*}\eta_j^{[\ell]}\bigg]\.\al\.\Jc\notag\\
   &\hspace*{20mm}=-(\im\la)^{-1}\.\Jc\al^*
     \int_{-\infty}^{\infty} \im(t-\la)^{-1}\.\d\tau_{\al,\be}(t)\.\al\.\Jc,\label{E:sum.RS.integral}
 \end{align}
from which, upon multiplying both sides of~\eqref{E:Im.Mla.estimate} by $\al\.\Jc$ from the right and by 
$-\Jc\al^*=(\al\Jc)^*$ from the left, we get the estimate
 \begin{equation}\label{E:EE14f}
  \int_{-\infty}^\infty \frac{\im(t-\la)^{-1}}{\im\la}\.\d\tau_{\al,\be}(t)\leq \frac{\im M_{N+1}(\la)}{\im\la},
  \quad \la\in\Cbb\stm\Rbb.
 \end{equation} 
Actually, this basic relation can be significantly improved if we employ Corollary~\ref{C:expansion} as it is done in 
the following theorem. This result together with the holomorphic property of $M_{N+1}(\la)$ provides an explicit form 
of an integral representation of $M_{N+1}(\la)$, which fully agrees with its Nevanlinna property, cf. 
\cite[Equality~(5.14)]{fG.erT00} and also \cite[Section~59]{niA.imG93}, \cite[Theorem~1.6]{bS19}.

\begin{theorem}\label{T:Im.Mla.integral}
 Let $\al,\be\in\Ga$ be given and Hypothesis~\ref{H:WAC} hold. Then the imaginary part of the Weyl--Titchmarsh 
 $M(\la)$-function admits the integral representation
  \begin{equation}\label{E:Im.Mla.integral}
   \im M_{N+1}(\la)=\im(\la)\.M^{[1]}+\int_{-\infty}^\infty \im(t-\la)^{-1}\.\d\tau_{\al,\be}(t)
   \qtext{for all $\la\in\Cbb\stm\Rbb$,}
  \end{equation}
 where $M^{[1]}\coloneq \lim_{\mu\to\infty} M_{N+1}(i\mu)/(i\mu)\geq0$. Consequently,  
  \begin{equation}\label{E:Mla.integral}
   M_{N+1}(\la)=M^{[0]}+\la\.M^{[1]}
                       +\int_{-\infty}^\infty\bigg(\frac{1}{t-\la}-\frac{t}{1+t^2}\bigg)\.\d\tau_{\al,\be}(t)
   \qtext{for all $\la\in\Cbb\stm\Rbb$,}
  \end{equation}
 where $M^{[0]}\coloneq \re M_{N+1}(i)$.
\end{theorem}
\begin{proof}
 Let $\la\in\Cbb\stm\Rbb$, $\nu\in\Cbb\stm(\Rbb\cup\{\la,\bla\})$, and $\xi\in\Cbb^n$ be arbitrary. Then one can easily 
 verify that the sequence $g\coloneq[\Xc(\la)-\Xc(\nu)]\.\xi\in\Cbb(\Izp)^{2n}$ is a solution of 
 problem~\eqref{E:EE15a} with $f\coloneq \la\.\Xc(\la)-\nu\.\Xc(\nu)$. Thus, by Corollary~\ref{C:expansion}, we have
  \begin{equation}\label{E:Im.Mla.intgr.*}
   \sum_{s\in\Iz} G_{k,s}(\bla)\.\Ps_s\.\big[\Xc_s(\la)-\Xc_s(\nu)\big]\.\xi
    =\sum_{j=1}^{r}\sum_{\ell=1}^{r_j} (\la_j-\bla)^{-1}\. d_j^{[\ell]}\.z_k^{[\ell]}(\la_j) 
  \end{equation}
 with the coefficients
  \begin{align*}
   d_j^{[\ell]}&=\sum_{k\in\Iz} z_k^{[\ell]*}(\la_j)\.\Ps_k\.\big[\Xc_k(\la)-\Xc_k(\nu)\big]\.\xi\\
    &\hspace*{-1.1mm}\overset{\eqref{E:lagrange}}{=}
      \frac{1}{\la_j-\la}\.\Big[\zeta^*\.Q^{[N+1]*}\.\Jc\.\Xc_{N+1}(\la)+\eta_j^{[\ell]*}\al\.\Xc_0(\la)\Big]\.\xi\\
    &\hspace{30mm}
      -\frac{1}{\la_j-\nu}\.\Big[\zeta^*\.Q^{[N+1]*}\.\Jc\.\Xc_{N+1}(\nu)+\eta_j^{[\ell]*}\al\.\Xc_0(\nu)\Big]\.\xi\\
    &=\frac{1}{\la_j-\la}\.\Big[-\zeta^*\.Q^{[N+1]*}\.\be^*\.\Xi(\la)+\eta_j^{[\ell]*}\Big]\.\xi\\
      &\hspace*{40mm}-\frac{1}{\la_j-\nu}\.\Big[-\zeta^*\.Q^{[N+1]*}\.\be^*\.\Xi(\nu)+\eta_j^{[\ell]*}\Big]\.\xi\\
    &=\Big(\frac{1}{\la_j-\la}-\frac{1}{\la_j-\nu}\Big)\.\eta_j^{[\ell]*}\.\xi,
  \end{align*}
 where $\zeta\in\Cbb^{2n}$ is such that $z_{N+1}^{[\ell]*}(\la_j)=Q^{[N+1]}\.\zeta$ analogously  
 to~\eqref{E:boundary.conditions.equiv}, while the $n\times n$ matrices $\Xi(\la),\Xi(\nu)$ are the same as in the 
 discussion concerning equality~\eqref{E:Mla.Xla.nu}. Especially, if we take $k=0$, then 
 from~\eqref{E:Im.Mla.intgr.*},  \eqref{E:EE5b}, and~\eqref{E:Mla.Xla.nu} we obtain
  \begin{align*}
   &-\Jc\.\al^*\bigg\{\frac{\im M_{N+1}(\la)}{\im\la}
                           -\frac{1}{\bla-\nu}\big[M_{N+1}^*(\la)-M_{N+1}(\nu)\big]\bigg\}\.\xi\\
    &\hspace*{30mm}
     =-\Jc\.\al^*\sum_{j=1}^{r}\sum_{\ell=1}^{r_j}
         \bigg[\frac{1}{\abs{\la_j-\la}^2}-\frac{1}{(\la_j-\bla)\.(\la_j-\nu)}\bigg]\.
                                                               \eta_j^{[\ell]}\.\eta_j^{[\ell]*}\.\xi.
  \end{align*}
 Since this equality is satisfied for all $\xi\in\Cbb^n$ it follows, similarly as in~\eqref{E:sum.RS.integral}, that
  \begin{align*}
   &\frac{\im M_{N+1}(\la)}{\im\la}-\frac{1}{\im \la}\int_{-\infty}^\infty \im(t-\la)^{-1}\.\d\tau_{\al,\be}(t)\\
   &\hspace*{20mm}
    =\frac{1}{\bla-\nu}\big[M_{N+1}^*(\la)-M_{N+1}(\nu)\big]
        -\sum_{j=1}^{r}\sum_{\ell=1}^{r_j}\frac{1}{(\la_j-\bla)\.(\la_j-\nu)}\.\eta_j^{[\ell]}\.\eta_j^{[\ell]*}.
  \end{align*}
 However, the left hand side does not depend on $\nu$, so it has to be true for any $\nu=i\.\mu$ such that 
 $\mu>\abs{\la}$. Therefore, as $\mu\to\infty$ the right-hand side tends to $M^{[1]}$, while the left-hand side 
 is positive semidefinite by~\eqref{E:EE14f}. This justifies equality~\eqref{E:Im.Mla.integral}. Consequently, since
 $M_{N+1}(\la)$ is holomorphic for all $\la\in\Cbb\stm\Rbb$, it is uniquely determined by its imaginary part up to 
 an additive Hermitian matrix, so it holds
  \begin{equation*}
   M_{N+1}(\la)=\la\.M^{[1]}+\int_{-\infty}^\infty \Big(\frac{1}{t-\la}+C\Big)\.\d\tau_{\al,\be}(t)
  \end{equation*}
 for some $C\in\Rbb$, which for the choice $\la=i$ leads to
 $\int_{-\infty}^\infty C\.\d\tau_{\al,\be}(t)=\re M_{N+1}(i)
 -\int_{-\infty}^\infty \frac{t}{1+t^2}\.\.\d\tau_{\al,\be}(t)$ implying equality~\eqref{E:Mla.integral}.
\end{proof}

\begin{example}
 For illustrative purposes, let us return to the system from the first part of Example~\ref{Ex:expansion}(i) with 
 $\al=(1\,\,\,\,0)=\be$. In that case we have the $M(\la)$-function
  \begin{equation*}
   M_{N+1}(\la)=-\frac{1}{\la\.v_{N+1}}, \quad\la\in\Cbb\stm\{0\},
  \end{equation*}
 and the corresponding spectral function
  \begin{equation*}
   \tau_{\al,\be}(t)
    =\begin{cases}
      0, & t\geq0,\\
      -1/v_{N+1}, & t<0.
     \end{cases}
  \end{equation*}
 Furthermore, $M^{[0]}=\re M_{N+1}(i)=0$ and $M^{[1]}=\lim_{\mu\to\infty} \frac{1}{\mu^2\.v_{N+1}}=0$, which together 
 with the equalities
  \begin{equation*}
   \int_{-\infty}^\infty \im(t-\la)^{-1}\.\d\tau_{\al,\be}(t)=-\frac{\im 1/\la}{v_{N+1}}
    =\frac{\im\la}{\abs{\la}^2\.v_{N+1}}
  \end{equation*}
 and  
  \begin{equation*}
   \int_{-\infty}^\infty \Big(\frac{1}{t-\la}-\frac{t}{1+t^2}\Big)\.\d\tau_{\al,\be}(t)=-\frac{1}{\la\.v_{N+1}}
  \end{equation*}
 yield the same integral representations of $\im M_{N+1}(\la)$ and $M_{N+1}(\la)$ as 
 displayed in~\eqref{E:Im.Mla.integral} and~\eqref{E:Mla.integral}.
\end{example}

Finally, we focus on the extension of the finite interval $\Iz=\onZ$ to the case $\oinftyZ$. Hence, let 
$N_0\in\oinftyZ$ be such that Hypothesis~\ref{H:WAC} is satisfied with $\Iz=\oinftyZ$ and $N$ replaced by $N_0$ in 
inequality~\eqref{E:Atk.ineq}. From the analysis of the so-called {\it Weyl disks} $D_{N}(\la)$ we know that 
$M_{N+1}(\la)$ lies on its boundary, which is known as the {\it Weyl circle}, see \cite[Section~3]{rSH.pZ14:JDEA}.
These disks are nested as $N$ tends to infinity, hence $M_{N+1}(\la)\in D_{N_0+1}(\la)$ for any 
$N\in[N_0,\infty)_\sZbb$ and, consequently, the sequence $\{M_{N+1}(\la)\}_{N\in[N_0,\infty)_\sZbb}$ is uniformly 
bounded. In particular, if we choose $\la=i$, then there exists a~matrix $T_0\in\Cbb^{n\times n}$ such that 
$\im M_{N+1}(i)\leq T_0$, which yields
 \begin{equation*}
  0\leq \int_{-\infty}^\infty (1+t^2)^{-1} \.\d\tau_{\al,\be}(t)
   =\int_{-\infty}^\infty \im(t-i)^{-1}\.\d\tau_{\al,\be}(t)\leq T_0
 \end{equation*}
by~\eqref{E:EE14f}. Consequently, for any $u>0$ we have
 \begin{align*}
  (1+u^2)^{-1}\d\tau_{\al,\be}(u)
   &=(1+u^2)^{-1}\.\int_0^{u}\.\d\tau_{\al,\be}(t)
    \leq \int_0^{u}\.(1+t^2)^{-1}\.\d\tau_{\al,\be}(t)\\
   &\leq \int_{-\infty}^{\infty}\.(1+t^2)^{-1}\.\d\tau_{\al,\be}(t)\leq T_0,
 \end{align*}
i.e., it holds
 \begin{equation}\label{E:EE14g}
  \tau_{\al,\be}(u)\leq (1+u^2)\.T_0,\quad u>0.
 \end{equation}
Similarly, for $u<0$ we get 
 \begin{equation*}
  -(1+u^2)^{-1}\d\tau_{\al,\be}(u)
   =(1+u^2)^{-1}\.\int_{u}^0\.\d\tau_{\al,\be}(t)
   \leq \int_{-\infty}^{\infty}\.(1+t^2)^{-1}\.\d\tau_{\al,\be}(t)\leq T_0,
 \end{equation*}
i.e., it holds
 \begin{equation}\label{E:EE14h}
  -\tau_{\al,\be}(u)\leq (1+u^2)\.T_0, \quad u<0.
 \end{equation}
Upon combining~\eqref{E:EE14g} and~\eqref{E:EE14h} we get the following statement.
 
\begin{theorem}\label{T:T7B}
 Let $\Iz=\oinftyZ$, matrices $\al,\be\in\Ga$ be given, and Hypothesis~\ref{H:WAC} hold with $N$ replaced by 
 some $N_0\in\oinftyZ$. Then there is a number $c>0$ independent of $\be$ and $N\in[N_0,\infty)_\sZbb$ such that
  \begin{equation}\label{E:EE14i}
   \abs{\tr \tau_{\al,\be}(t)}\leq c\.(1+t^2) \qtext{for any $t\in\Rbb$.}
  \end{equation}
\end{theorem}

Inequality~\eqref{E:EE14i} yields a bound for the absolute value of every diagonal element of $\tau_{\al,\be}(t)$. 
Actually, it provides a bound of all elements of $\tau_{\al,\be}(t)$ because of the Hermitian property of 
this matrix-valued function. Consequently, if we take an arbitrary compact subinterval $[a,b]\subseteq\Rbb$, then by 
using similar arguments as in~\cite[Lemma~2 on p.~221]{ipN55} it can be shown that for $\la\in[a,b]$ there is a 
sequence $N_0<N_1<N_2<\dots$ with $N_j\to\infty$ and a~nondecreasing right-continuous Hermitian matrix-valued function 
$\tau(\cdot)$ such that
 \begin{equation*}
  \tau_{\al,\be,N_j}(\la)\to\tau(\la) \qtext{as $j\to\infty$}
 \end{equation*}
whenever $\tau(\cdot)$ is continuous at $\la$, i.e., except for at most a countable subset of $\Rbb$. Furthermore, this 
{\it limiting spectral function} also satisfies the inequalities
 \begin{equation*}
  \sgn(\la)\.\tau(\la)\geq0 \qtextq{and}
  \abs{\tr\tau(\la)}\leq c\.(1+\la^2) \qtext{for all $\la\in\Rbb$.}
 \end{equation*} 

The limiting spectral function shall be utilized in our subsequent research concerning an investigation of the 
spectrum of self-adjoint linear relations associated with the mapping $\mL(\cdot)$. At this moment we show only a 
certain analogue of Parserval's identity~\eqref{E:parseval}, in which we connect $\tau(\cdot)$ and 
the norm of compactly supported solutions of~\Slaf{0}{f}. It is based on the previous results of this section, in which 
we take $N\to\infty$, compare with \cite[Theorem~9.9.2]{fvA64}. Actually, these solutions form the domain of 
a~preminimal relation associated with $\mL(\cdot)$, see \cite[Section~5]{slC.pZ15} and \cite{pZ.slC16} for more 
details.

\begin{theorem}\label{T:T7C}
 Let $\Iz=\oinftyZ$, a matrix $\al\in\Ga$ be given and $\hz\in\Cbb(\Izp)^{2n}$ be such that $\hz_0=0=\hz_k$ for all 
 $k\in\Iz$ large enough. If Hypothesis~\ref{H:WAC} holds for $N$ replaced by some $N_0\in\Iz$ and there is $f\in\ltp$ 
 such that $\hz$ solves system~\Slaf{0}{f} on $\Iz$, then for the $n$-vector valued function $\om(\la)$ given as 
 the~finite sum 
  \begin{equation*}
   \om(\la)\coloneq \sum_{k\in\Iz} \tZ_k^*(\la)\.\Ps_k\.\hz_k
  \end{equation*}
 we have 
  \begin{equation*}
   \normP{\hz}^2=\int_{-\infty}^\infty \om^*(\la)\.\d\tau(\la)\.\om(\la).
  \end{equation*}
\end{theorem}
\begin{proof}
 Let $a>0$ be arbitrary and $N\in\Iz$ be such that Hypothesis~\ref{H:WAC} holds on the finite discrete interval $\onZ$ 
 and $\hz_k=0$ for all $k\in[N+1,\infty)_\sZbb$. Then $\Ps_k\.f_k\equiv0$ on $[N+1,\infty)_\sZbb$ and these 
 two sequences $\hz$ and $f$ satisfy the assumptions of Lemma~\ref{T:T9} for any $\be\in\Ga$ with the finite discrete 
 interval corresponding to $\onZ$. Hence, from~\eqref{E:EE25} we get 
  \begin{equation}\label{E:EE26}
   \normP{\hz}^2-\sum_{\abs{\la_j}\leq a}\sum_{\ell=1}^{r_j} \abs[\big]{c_j^{[\ell]}}^2\leq a^{-2}\.\normP{f}^2.
  \end{equation}
 At the same time, from the formula for $c_j^{[\ell]}$ displayed in~\eqref{E:expansion.coeff} we obtain
  \begin{equation*}
   c_j^{[\ell]}=\eta_j^{[\ell]*}\.\om(\la_j)
  \end{equation*}
 for $\eta_j^{[\ell]}\in\ker\be\.\tZ_{N+1}(\la_j)$ being the vectors from the orthonormalization process discussed 
 before Theorem~\ref{T:orthonormal.set}. So, if $\pm a$ are not points of discontinuity of $\tau_{\al,\be,N}$, i.e., 
 they are not eigenvalues of~\eqref{E:Sla} on $\onZ$, it follows 
  \begin{align}\label{E:EE27}
   \sum_{\abs{\la_j}\leq a}\sum_{\ell=1}^{r_j} \abs[\big]{c_j^{[\ell]}}^2
    &=\sum_{\abs{\la_j}\leq a}\sum_{\ell=1}^{r_j} \om^*(\la_j)\.\eta_j^{[\ell]}\.\eta_j^{[\ell]*}\.\om(\la_j)\notag\\
    &=\int_{-a}^a \om^*(\la)\.\d\tau_{\al,\be,N}(\la)\.\om(\la).
  \end{align}
 If $\pm a$ are not eigenvalues of system~\eqref{E:Sla} on $[0,N_j]_\sZbb$ for the sequence of endpoints used
 above for the existence of the limiting spectral function $\tau$, i.e., they are not points of discontinuity of 
 the corresponding $\tau_{\al,\be,N_j}$ and $\tau$, then we can let $N\to\infty$ in~\eqref{E:EE26}--\eqref{E:EE27}, 
 which yields
  \begin{equation*}
   \normP{\hz}^2-\int_{-a}^a \om^*(\la)\.\d\tau(\la)\.\om(\la)\leq a^{-2}\.\normP{f}^2
  \end{equation*}
 and, consequently, as $a$ tends to infinity we achieve the desired conclusion.
\end{proof}


\section*{Acknowledgments}
\addcontentsline{toc}{section}{Acknowledgments}
The research was supported by the Czech Science Foundation under Grant GA19-01246S. The author is grateful to the 
anonymous referee for detailed reading of the manuscript and her/his comments.



\end{document}